\documentclass[a4paper]{jpconf}

\usepackage{lmodern}
\usepackage{amsmath}
\usepackage{amssymb}
\usepackage{graphicx}
\usepackage{nomencl}
\usepackage{cite}
\usepackage{esint}
\usepackage{geometry}
\begin{document}

\title{Semi-analytical calculation of the singular and hypersingular integrals
for discrete Helmholtz operators in 2D BEM}

\author{Andrea Cagliero$^1$}
\address{$^1$ The author carried out this study while affiliated with the Microwaves Department of IMT Atlantique, Institut Mines-T\'{e}l\'{e}com, and with the Laboratory for Science and Technologies of Information, Communication and Knowledge Lab-STICC (CNRS), Brest, F-29238, France.}

\ead{andrea.cagliero@edu.unito.it}

\begin{abstract}
Approximate solutions to elliptic partial differential equations with
known kernel can be obtained via the boundary element method (BEM)
by discretizing the corresponding boundary integral operators and
solving the resulting linear system of algebraic equations. Due
to the presence of singular and hypersingular integrals, the evaluation
of the operator matrix entries requires the use of regularization
techniques. In this work, the singular and hypersingular integrals
associated with first-order Galerkin discrete boundary operators for
the two-dimensional Helmholtz equation are reduced to quasi-closed-form
expressions. The obtained formulas may prove useful for the implementation
of the BEM in two-dimensional electromagnetic, acoustic and quantum
mechanical problems.
\end{abstract}

\section{Introduction}

The boundary element method (BEM) is one of the most common numerical
techniques to solve elliptic boundary value problems \cite{Steinbach2008,Sauter2011}.
Making use of the kernel of the considered partial differential equation,
integral theorems are employed to express the solution in terms of
bounded operators on Sobolev spaces. The Cauchy data of the problem
is then obtained by applying an appropriate discretization scheme,
which consists in projecting the solution onto finite dimensional
trial spaces, and by numerically solving the resulting algebraic equations.
Among the possible discretization strategies, the weak formulation
known as the symmetric Galerkin method has been largely considered
in the literature (see \cite{Sutradhar2008} and references therein). 

In contrast to the finite element method (FEM), the BEM has the advantage
of only requiring the discretization of the boundary of the physical
domain without the need to introduce any truncation in open-region
problems. However, since most boundary integral operators are singular,
regularization procedures must be taken into account. In the present
study, a semi-analytical approach is proposed to evaluate all the
possible singular integrals arising from the first-order Galerkin
discretization of the boundary operators for the two-dimensional Helmholtz
equation. These results may be relevant for BEM applications in electromagnetism
and acoustics \cite{Nedelec2001,Jin2015} as well as in quantum mechanics
\cite{Ramdas2002}. 

The manuscript is organized as follows. In Section \ref{sec:problemStatement},
the boundary integral operators are introduced in both their continuous
and discrete forms; the two-dimensional Helmholtz kernel derivatives
and the linear basis functions are consequently defined. Section \ref{sec:singleLayer}
is concerned with the calculation of the singular integrals occurring
in the discrete single layer operator. In Section \ref{sec:doubleLayer}
and \ref{sec:directMethod}, the same analysis is carried out for
the discrete double layer operators and for the discrete hypersingular
operator, respectively. An alternative approach based on the variational
formulation for the hypersingular operator is reported in Section
\ref{sec:variationalApproach}. Finally, Section \ref{sec:conclusions}
is devoted to the conclusions.

\section{Problem statement\label{sec:problemStatement}}

Let $g\left(\mathbf{r},\mathbf{r}'\right)$ be the kernel
of an elliptic partial differential equation over the domain $\varOmega$
and $f\left(\mathbf{r}\right)$ a well-behaved function defined
on $S\equiv\partial\varOmega$. The four boundary integral operators
known as single layer, double layer, adjoint double layer and hypersingular
are defined, respectively, as follows \cite{Nedelec2001,Steinbach2008}:

\begin{align}
\hat{S}\left[f\right]\left(\mathbf{r}\right) & \equiv\fint_{S}d\mathbf{r}'g\left(\mathbf{r},\mathbf{r}'\right)f\left(\mathbf{r}'\right);\label{eq:singleLayer}\\
\hat{D}\left[f\right]\left(\mathbf{r}\right) & \equiv\fint_{S}d\mathbf{r}'\frac{\partial g\left(\mathbf{r},\mathbf{r}'\right)}{\partial n'}f\left(\mathbf{r}'\right)=\fint_{S}d\mathbf{r}'\nabla'g\left(\mathbf{r},\mathbf{r}'\right)\cdot\mathbf{n}'\,f\left(\mathbf{r}'\right);\\
\hat{D}^{\dagger}\left[f\right]\left(\mathbf{r}\right) & \equiv\fint_{S}d\mathbf{r}'\frac{\partial g\left(\mathbf{r},\mathbf{r}'\right)}{\partial n}f\left(\mathbf{r}'\right)=\fint_{S}d\mathbf{r}'\nabla g\left(\mathbf{r},\mathbf{r}'\right)\cdot\mathbf{n}\,f\left(\mathbf{r}'\right);\\
\hat{N}\left[f\right]\left(\mathbf{r}\right) & \equiv\int_{S}d\mathbf{r}'\frac{\partial^{2}g\left(\mathbf{r},\mathbf{r}'\right)}{\partial n\partial n'}f\left(\mathbf{r}'\right)=\int_{S}d\mathbf{r}'\nabla\left[\nabla'g\left(\mathbf{r},\mathbf{r}'\right)\cdot\mathbf{n}'\right]\cdot\mathbf{n}\,f\left(\mathbf{r}'\right),\label{eq:hyperSingular}
\end{align}
where the symbol $\fint$ stands for the Cauchy principal value integral
and $\mathbf{n}$, $\mathbf{n}'$ are the outward pointing
unit normals to $S$ evaluated at $\mathbf{r}$ and $\mathbf{r}'$,
respectively. A formal solution to the considered elliptic partial
differential equation is often provided in terms of (\ref{eq:singleLayer})-(\ref{eq:hyperSingular}).
In order to numerically implement these integral operators, a common
strategy is to discretize the surface $S$ into a collection of simplices
$\left\{ S_{n}\right\} $. The unknown solution $f\left(\mathbf{r}\right)$
is then expanded over a set of basis functions defined on specific
groups of neighbor simplices. For instance, when the basis functions
are a given set of interpolation polynomials $\left\{ p_{j}\left(\mathbf{r}\right)\right\} $,
we have: 
\begin{equation}
f\left(\mathbf{r}\right)=\sum_{j}\alpha_{j}p_{j}\left(\mathbf{r}\right),
\end{equation}
with $\alpha_{j}$ representing the value of the function $f\left(\mathbf{r}\right)$
at the $j$-th mesh node. The $j$-th basis function $p_{j}\left(\mathbf{r}\right)$
is defined on the set of simplices $\left\{ S_{n}\right\} $ that
share the $j$-th mesh node, hereinafter referred to as $\left\{ n\in j\right\} $,
and vanishes out of its defining domain, so that:
\begin{equation}
\int_{S}d\mathbf{r}'p_{j}\left(\mathbf{r}'\right)=\sum_{n\in j}\int_{S_{n}}d\mathbf{r}'p_{j}^{n}\left(\mathbf{r}'\right),
\end{equation}
where $p_{j}^{n}\left(\mathbf{r}'\right)$ is the restriction
of the $j$-th basis function to the $n$-th simplex. According to
the well-known Galerkin approach \cite{Gibson2015,Jin2015}, the same
set of basis functions may be used to symmetrize the discrete version
of the surface integral operators, leading to:
\begin{align}
S_{ij} & \equiv\sum_{m\in i}\,\sum_{n\in j}\int_{S_{m}}d\mathbf{r}\fint_{S_{n}}d\mathbf{r}'g\left(\mathbf{r},\mathbf{r}'\right)p_{i}^{m}\left(\mathbf{r}\right)p_{j}^{n}\left(\mathbf{r}'\right);\label{eq:singleLayer_Matrix}\\
D_{ij} & \equiv\sum_{m\in i}\,\sum_{n\in j}\int_{S_{m}}d\mathbf{r}\fint_{S_{n}}d\mathbf{r}'\frac{\partial g\left(\mathbf{r},\mathbf{r}'\right)}{\partial n'}p_{i}^{m}\left(\mathbf{r}\right)p_{j}^{n}\left(\mathbf{r}'\right);\label{eq:doubleLayer_Matrix}\\
D_{ij}^{\dagger} & \equiv\sum_{m\in i}\,\sum_{n\in j}\int_{S_{m}}d\mathbf{r}\fint_{S_{n}}d\mathbf{r}'\frac{\partial g\left(\mathbf{r},\mathbf{r}'\right)}{\partial n}p_{i}^{m}\left(\mathbf{r}\right)p_{j}^{n}\left(\mathbf{r}'\right);\label{eq:adjointDoubleLayer_Matrix}\\
N_{ij} & \equiv\sum_{m\in i}\,\sum_{n\in j}\int_{S_{m}}d\mathbf{r}\int_{S_{n}}d\mathbf{r}'\frac{\partial^{2}g\left(\mathbf{r},\mathbf{r}'\right)}{\partial n\partial n'}p_{i}^{m}\left(\mathbf{r}\right)p_{j}^{n}\left(\mathbf{r}'\right).\label{eq:hyperSingular_Matrix}
\end{align}

\begin{figure}
\begin{centering}
\includegraphics[width=1\textwidth]{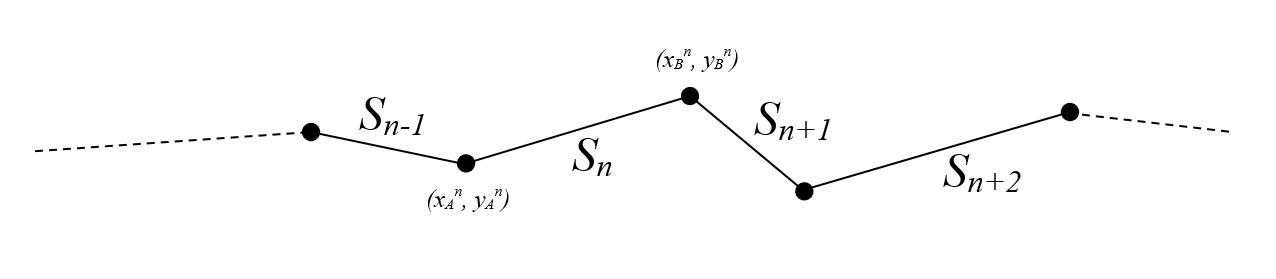}
\par\end{centering}
\caption{Collection of segments $\left\{ S_{n}\right\} $ resulting from the
discretization of a 2D curve.\label{fig:segments}}

\end{figure}

Let us now focus on the Helmholtz equation over the 2D region $\varOmega\subset\mathbb{R}^{2}$
whose boundary $S$ is a piecewise smooth closed curve. The kernel
of the equation is expressed in terms of the Hankel function \cite{Abramowitz1972}:
\begin{equation}
g\left(x,y;x',y'\right)=\frac{i}{4}H_{0}^{(1,2)}\left[k\sqrt{\left(x-x'\right)^{2}+\left(y-y'\right)^{2}}\right]\label{eq:Green}
\end{equation}
and its normal derivatives are given by:
\begin{align}
\frac{\partial g\left(x,y;x',y'\right)}{\partial n'} & =\frac{ik}{4R}H_{1}^{(1,2)}\left(kR\right)\left(\mathbf{R}\cdot\mathbf{n}'\right);\label{eq:Dnjg}\\
\frac{\partial g\left(x,y;x',y'\right)}{\partial n} & =-\frac{ik}{4R}H_{1}^{(1,2)}\left(kR\right)\left(\mathbf{R}\cdot\mathbf{n}\right);\label{eq:Dnig}\\
\frac{\partial^{2}g\left(x,y;x',y'\right)}{\partial n\partial n'} & =\frac{ik}{4R^{2}}\left[R\,H_{1}^{(1,2)}\left(kR\right)\left(\mathbf{n}\cdot\mathbf{n}'\right)-kH_{2}^{(1,2)}\left(kR\right)\left(\mathbf{R}\cdot\mathbf{n}\right)\left(\mathbf{R}\cdot\mathbf{n}'\right)\right],\label{eq:DniDnjg}
\end{align}
with $\mathbf{R}\equiv\mathbf{r}-\mathbf{r}'=\left(x-x',y-y'\right)$
and $k$ representing the wave number. The boundary curve $S$ can
be discretized into a collection of segments $\left\{ S_{n}\right\} $
with lengths $\left\{ l_{n}\right\} $ and extrema $\left\{ \mathbf{r}_{A}^{n};\:\mathbf{r}_{B}^{n}\right\} =\left\{ \left(x_{A}^{n},y_{A}^{n}\right);\:\left(x_{B}^{n},y_{B}^{n}\right)\right\} $,
as depicted in Figure \ref{fig:segments}.
 First order basis functions (triangular functions) are defined over
pairs of adjacent segments and vary linearly from zero at the outer
extrema to unity at the common vertex \cite{Gibson2015,Jin2015}.
By introducing the local variable $t_{n}\in\left[0,1\right]$, which
makes it possible to represent an arbitrary point $\mathbf{r}=\left(x,y\right)$
on the $n$-th segment in parametric form: 
\begin{equation}
\mathbf{r}\left(t_{n}\right)=\mathbf{r}_{A}^{n}+\left(\mathbf{r}_{B}^{n}-\mathbf{r}_{A}^{n}\right)t_{n},\label{eq:parametricCoordinates}
\end{equation}
the restriction of the $j$-th triangular basis function to the $n$-th
segment can be written as follows:%

\begin{equation}
p_{j}^{n}\left[\mathbf{r}\left(t_{n}\right)\right]\equiv\begin{cases}
1-t_{n} & \mathrm{if}\:\mathbf{r}_{j}=\mathbf{r}_{A}^{n};\\
t_{n} & \mathrm{if}\:\mathbf{r}_{j}=\mathbf{r}_{B}^{n},
\end{cases}\label{eq:basisFunctions}
\end{equation}
where $\mathbf{r}_{j}=\left(x_{j},y_{j}\right)$ identifies the
coordinates of the $j$-th mesh node.

\begin{figure}
\begin{centering}
\includegraphics[width=1\textwidth]{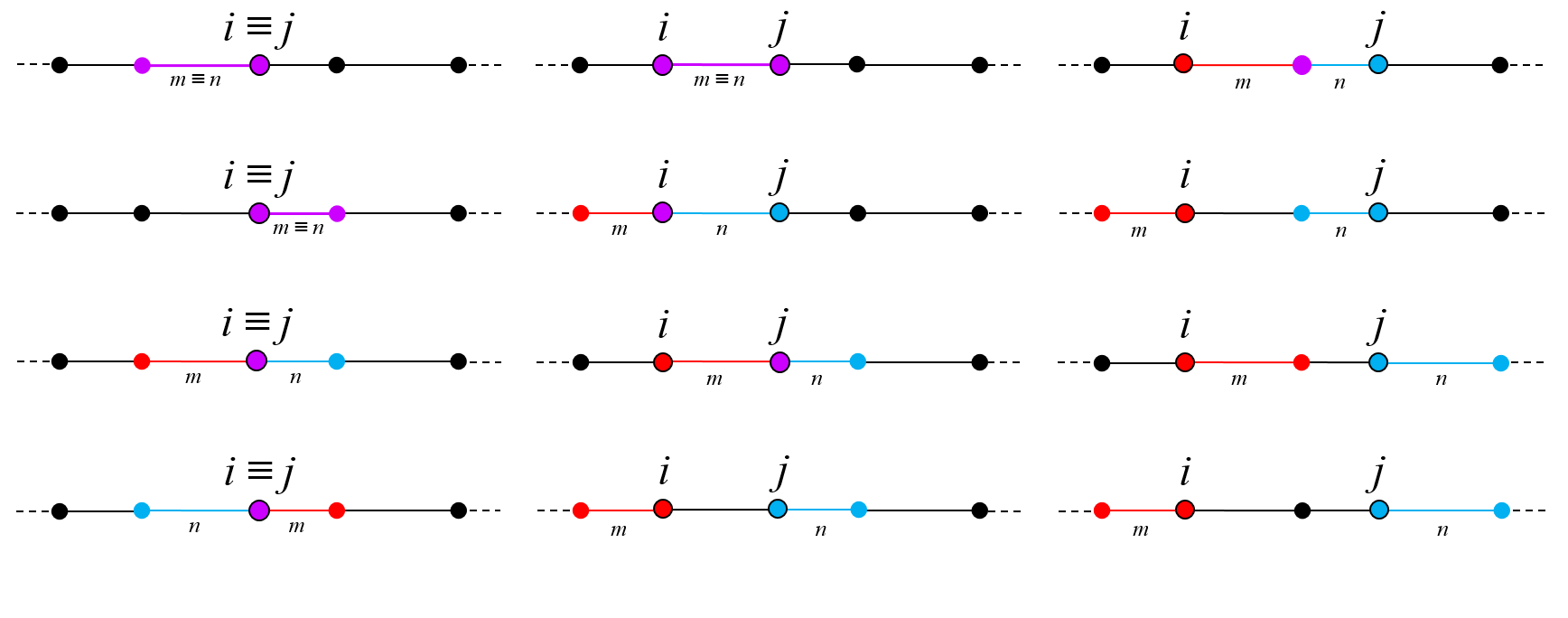}
\par\end{centering}
\caption{Graphical representation of the four double integrations occurring
in each of (\ref{eq:singleLayer_Matrix})-(\ref{eq:hyperSingular_Matrix})
when the $i$-th and $j$-th mesh nodes coincide (left column), when
they are first neighbors (central column) and when they are second-neighbors
(right column). Angles between segments are not displayed here (the
reader is referred to Figure \ref{fig:segments}). \label{fig:combinations}}
\end{figure}

In the present scenario, each of the discrete operators $S_{ij}$,
$D_{ij}$, $D_{ij}^{\dagger}$ and $N_{ij}$ defined in (\ref{eq:singleLayer_Matrix})-(\ref{eq:hyperSingular_Matrix})
consists of a sum of four double integrals over the pairs of segments
$\left(S_{m},S_{n}\right)\in\left\{ m\in i\right\} \times\left\{ n\in j\right\} $.
A graphical representation of these four terms is reported in Figure
\ref{fig:combinations} for three different choices of the mesh nodes
$i$ and $j$.
 When $S_{m}$ and $S_{n}$ do not share any vertex, all such double
integrals can be computed by Gauss-Legendre quadrature rules \cite{Abramowitz1972}.
The remaining double integrals with $S_{m}=S_{n}$ in (\ref{eq:singleLayer_Matrix})
and (\ref{eq:hyperSingular_Matrix}) as well as some of the cases
where $S_{m}$ and $S_{n}$ share only one vertex (i.e., adjacent
segments) are singular and must be treated with due care. In the following,
analytical integration formulas will be applied to each of these singular
double integrals in order to recast them into a quasi-closed-form
expression involving Bessel-related functions and two regular single
integrals that can be easily solved numerically. 

\section{Single layer operator\label{sec:singleLayer}}

\subsection{Integration over coincident elements\label{subsec:coincidentIntSingleLayer}}

\begin{figure}[t]
\begin{centering}
\includegraphics[width=1\textwidth]{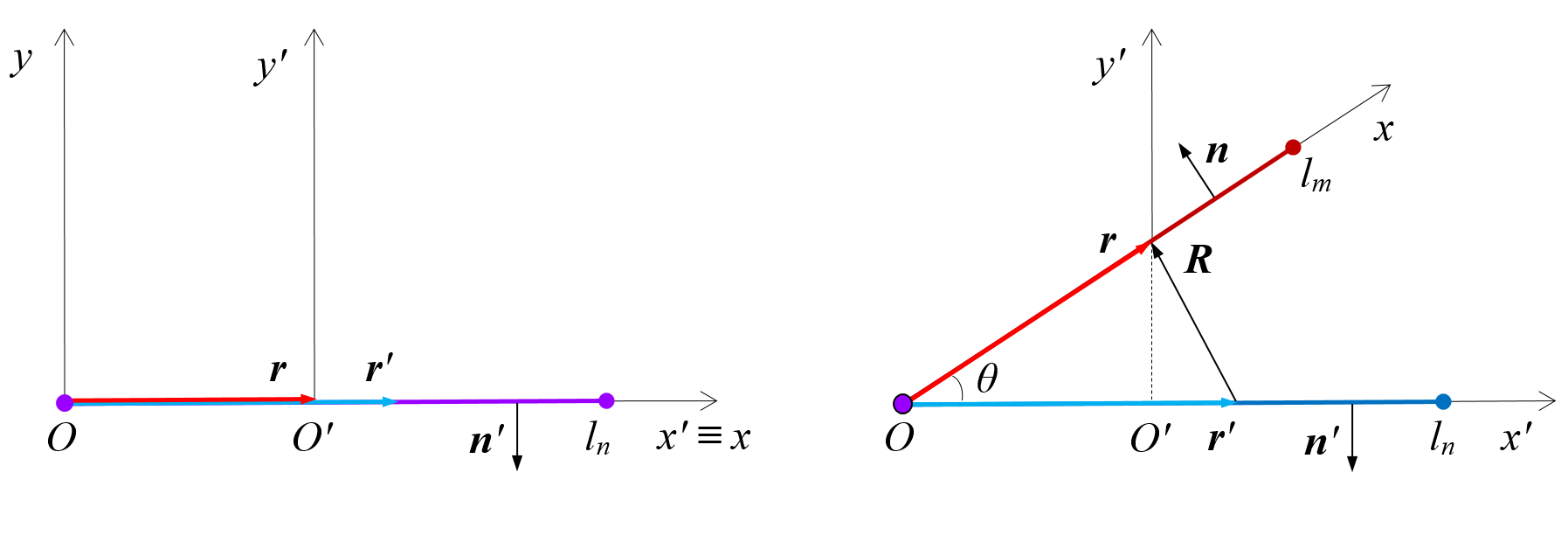}
\par\end{centering}
\caption{Sketches of the reference frames $Oxy$, $O'x'y'$ relative to the
coordinates $\mathbf{r}$ and $\mathbf{r}'$ for the double
integration over coincident (left) and adjacent (right) segments $S_{m}$
and $S_{n}$, respectively. \label{fig:referenceFrames}}
\end{figure}

When $S_{m}=S_{n}$, the singular double integrals in (\ref{eq:singleLayer_Matrix})
are of the form:
\begin{equation}
\frac{i}{4}\int_{S_{n}}d\mathbf{r}\fint_{S_{n}}d\mathbf{r}'H_{0}^{(1,2)}\left(k\left|\mathbf{r}-\mathbf{r}'\right|\right)p_{i}^{n}\left(\mathbf{r}\right)p_{j}^{n}\left(\mathbf{r}'\right),\label{eq:singularIntegral}
\end{equation}
where the linear functions $p_{i}^{n}\left(\mathbf{r}\right)$,
$p_{j}^{n}\left(\mathbf{r}'\right)$ are defined in (\ref{eq:basisFunctions})
via (\ref{eq:parametricCoordinates}). A useful choice for the reference
frames $Oxy$, $O'x'y'$ relative to the two integration coordinates
$\mathbf{r}$ and $\mathbf{r}'$ is to have both the $x$-axis
and the $x'$-axis lie along the segment $S_{n}$, with origins $O$,
$O'$ at $\left(x_{A}^{n},y_{A}^{n}\right)$ and $\mathbf{r}$,
respectively (see the left side of Figure \ref{fig:referenceFrames}).
With this convention, the $y$ and $y'$ coordinates become unnecessary
and the integral (\ref{eq:singularIntegral}) reduces to:
\begin{equation}
\frac{i}{4}\int_{0}^{l_{n}}dx\fint_{-x}^{l_{n}-x}dx'H_{0}^{(1,2)}\left(k\left|x'\right|\right)p_{i}^{n}\left(x\right)\widetilde{p}_{j}^{n}\left(x',x\right),\label{eq:singIntRewritten}
\end{equation}
where:
\begin{equation}
p_{i}^{n}\left(x\right)\equiv\begin{cases}
1-\frac{x}{l_{n}} & \mathrm{if}\:x_{i}=0;\\
\frac{x}{l_{n}} & \mathrm{if}\:x_{i}=l_{n},
\end{cases}\label{eq:bfi}
\end{equation}
\begin{equation}
\widetilde{p}_{j}^{n}\left(x',x\right)\equiv\begin{cases}
1-\frac{\left(x'+x\right)}{l_{n}} & \mathrm{if}\:x_{j}'=-x;\\
\frac{\left(x'+x\right)}{l_{n}} & \mathrm{if}\:x_{j}'=l_{n}-x.
\end{cases}\label{eq:bfj}
\end{equation}
Four kinds of integrals are obtained from (\ref{eq:singIntRewritten}),
(\ref{eq:bfi}) and (\ref{eq:bfj}), namely:
\begin{align}
I_{11}^{n} & \equiv\frac{i}{4}\int_{0}^{l_{n}}dx\fint_{-x}^{l_{n}-x}dx'H_{0}^{(1,2)}\left(k\left|x'\right|\right)\left(1-\frac{x}{l_{n}}\right)\left[1-\frac{\left(x'+x\right)}{l_{n}}\right];\label{eq:I11}\\
I_{12}^{n} & \equiv\frac{i}{4}\int_{0}^{l_{n}}dx\fint_{-x}^{l_{n}-x}dx'H_{0}^{(1,2)}\left(k\left|x'\right|\right)\left(1-\frac{x}{l_{n}}\right)\left[\frac{\left(x'+x\right)}{l_{n}}\right];
\end{align}
\begin{align}
I_{21}^{n} & \equiv\frac{i}{4}\int_{0}^{l_{n}}dx\fint_{-x}^{l_{n}-x}dx'H_{0}^{(1,2)}\left(k\left|x'\right|\right)\left(\frac{x}{l_{n}}\right)\left[1-\frac{\left(x'+x\right)}{l_{n}}\right];\\
I_{22}^{n} & \equiv\frac{i}{4}\int_{0}^{l_{n}}dx\fint_{-x}^{l_{n}-x}dx'H_{0}^{(1,2)}\left(k\left|x'\right|\right)\left(\frac{x}{l_{n}}\right)\left[\frac{\left(x'+x\right)}{l_{n}}\right].\label{eq:I22}
\end{align}
The integrals $I_{11}^{n}$ and $I_{22}^{n}$ may be linked to the
graphical representations in the second and first row of the left
column of Figure \ref{fig:combinations}, respectively. In other words,
both $I_{11}^{n}$ and $I_{22}^{n}$ will only arise when the mesh
nodes $i$ and $j$ coincide. Similarly, $I_{12}^{n}$ can be associated
to the first sketch in the central column of Figure \ref{fig:combinations},
as would be for $I_{21}^{n}$ if we exchanged the node labels: these
two integrals come into play when $i$ and $j$ are first neighbors. 

Making use of the following definitions:
\begin{align}
I_{\alpha}^{n}\left(x\right) & \equiv\fint_{-x}^{l_{n}-x}dx'H_{0}^{(1,2)}\left(k\left|x'\right|\right);\label{eq:Ialphan}\\
I_{\beta}^{n}\left(x\right) & \equiv\frac{x}{l_{n}}\fint_{-x}^{l_{n}-x}dx'H_{0}^{(1,2)}\left(k\left|x'\right|\right);\label{eq:Ibeta}\\
I_{\gamma}^{n}\left(x\right) & \equiv\fint_{-x}^{l_{n}-x}dx'H_{0}^{(1,2)}\left(k\left|x'\right|\right)\frac{\left(x'+x\right)}{l_{n}};\\
I_{\delta}^{n}\left(x\right) & \equiv\frac{x}{l_{n}}\fint_{-x}^{l_{n}-x}dx'H_{0}^{(1,2)}\left(k\left|x'\right|\right)\frac{\left(x'+x\right)}{l_{n}},\label{eq:Idelta}
\end{align}
the above integrals are rewritten as:
\begin{align}
I_{11}^{n} & =\frac{i}{4}\int_{0}^{l_{n}}dx\left[I_{\alpha}^{n}\left(x\right)-I_{\beta}^{n}\left(x\right)-I_{\gamma}^{n}\left(x\right)+I_{\delta}^{n}\left(x\right)\right];\label{eq:I11new}\\
I_{12}^{n} & =\frac{i}{4}\int_{0}^{l_{n}}dx\left[I_{\gamma}^{n}\left(x\right)-I_{\delta}^{n}\left(x\right)\right];\\
I_{21}^{n} & =\frac{i}{4}\int_{0}^{l_{n}}dx\left[I_{\beta}^{n}\left(x\right)-I_{\delta}^{n}\left(x\right)\right];\\
I_{22}^{n} & =\frac{i}{4}\int_{0}^{l_{n}}dx\,I_{\delta}^{n}\left(x\right).\label{eq:I22new}
\end{align}
Let us focus on $I_{\alpha}^{n}\left(x\right)$:
\begin{align}
I_{\alpha}^{n}\left(x\right) & =\int_{-x}^{0}dx'H_{0}^{(1,2)}\left(-kx'\right)+\int_{0}^{l_{n}-x}dx'H_{0}^{(1,2)}\left(kx'\right)\nonumber \\
 & =\frac{1}{k}\int_{0}^{kx}dx'H_{0}^{(1,2)}\left(x'\right)+\frac{1}{k}\int_{0}^{k\left(l_{n}-x\right)}dx'H_{0}^{(1,2)}\left(x'\right),
\end{align}
where the changes of variable $x'\rightarrow-kx'$ and $x'\rightarrow kx'$
have been employed to transform the first and second integrals, respectively.
Now, applying the same strategy to (\ref{eq:Ibeta})-(\ref{eq:Idelta})
and introducing the useful definitions:
\begin{align}
\mathcal{I}_{0}\left(\sigma\right) & \equiv\int_{0}^{\sigma}dx'H_{0}^{(1,2)}\left(x'\right);\label{eq:I0}\\
\mathcal{I}_{1}\left(\sigma\right) & \equiv\int_{0}^{\sigma}dx'H_{0}^{(1,2)}\left(x'\right)\,x',\label{eq:I1}
\end{align}
we get:
\begin{align}
I_{\alpha}^{n}\left(x\right) & =\frac{1}{k}\left\{ \mathcal{I}_{0}\left(kx\right)+\mathcal{I}_{0}\left[k\left(l_{n}-x\right)\right]\right\} ;\label{eq:Ialpha}\\
I_{\beta}^{n}\left(x\right) & =\frac{x}{l_{n}}I_{\alpha}^{n}\left(x\right);\label{eq:IbetaNew}\\
I_{\gamma}^{n}\left(x\right) & =I_{\beta}^{n}\left(x\right)+\frac{1}{k^{2}l_{n}}\left\{ \mathcal{I}_{1}\left[k\left(l_{n}-x\right)\right]-\mathcal{I}_{1}\left(kx\right)\right\} ;\label{eq:Igamma}\\
I_{\delta}^{n}\left(x\right) & =\frac{x}{l_{n}}I_{\gamma}^{n}\left(x\right).\label{eq:IdeltaNew}
\end{align}

\begin{figure}
\begin{centering}
\includegraphics[width=1\textwidth]{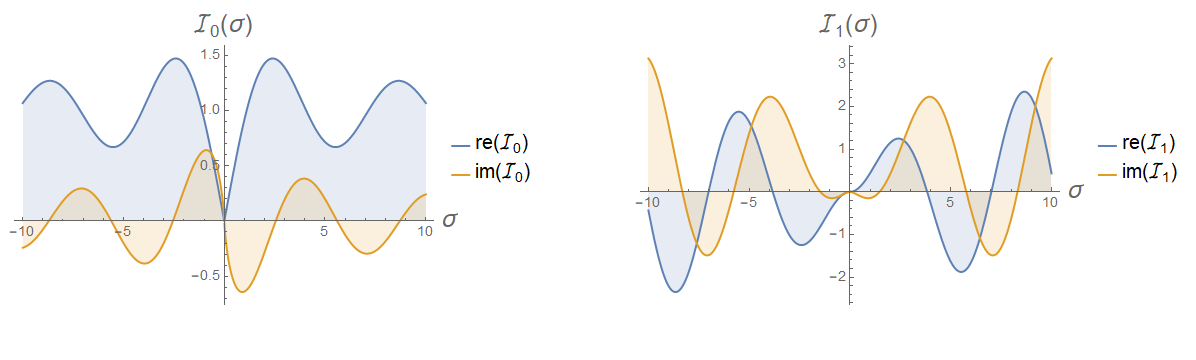}
\par\end{centering}
\caption{Plots of the functions $\mathcal{I}_{0}\left(\sigma\right)$ and $\mathcal{I}_{1}\left(\sigma\right)$
defined in (\ref{eq:I0}) and (\ref{eq:I1}) when only the first kind
Hankel function $H_{0}^{(1)}\left(x'\right)$ is considered.\label{fig:profiles}}
\end{figure}

In order to determine a closed-form expression for the integrals (\ref{eq:I0})
and (\ref{eq:I1}), reference is made to some tabulated formulas for
Bessel functions of the first and second kind \cite{Gradshteyn2007}:
\begin{align}
\int_{0}^{1}ds\,J_{\nu}\left(\sigma s\right)s^{\nu} & =2^{\nu-1}\sigma^{-\nu}\sqrt{\pi}\,\Gamma\left(\nu+\frac{1}{2}\right)\left[J_{\nu}\left(\sigma\right)\mathbf{H}_{\nu-1}\left(\sigma\right)-\mathbf{H}_{\nu}\left(\sigma\right)J_{\nu-1}\left(\sigma\right)\right];\\
\int_{0}^{1}ds\,Y_{\nu}\left(\sigma s\right)s^{\nu} & =2^{\nu-1}\sigma^{-\nu}\sqrt{\pi}\,\Gamma\left(\nu+\frac{1}{2}\right)\left[Y_{\nu}\left(\sigma\right)\mathbf{H}_{\nu-1}\left(\sigma\right)-\mathbf{H}_{\nu}\left(\sigma\right)Y_{\nu-1}\left(\sigma\right)\right];\\
\int_{0}^{1}ds\,J_{\nu}\left(\sigma s\right)s^{\nu+1} & =\sigma^{-1}J_{\nu+1}\left(\sigma\right);\\
\int_{0}^{1}ds\,Y_{\nu}\left(\sigma s\right)s^{\nu+1} & =\sigma^{-1}Y_{\nu+1}\left(\sigma\right)+2^{\nu+1}\sigma^{-\nu-2}\pi^{-1}\Gamma\left(\nu+1\right),
\end{align}
where $\mathbf{H}_{\nu}$ represents the Struve function of order
$\nu$ and $\Gamma$ is the Gamma function. If we combine the previous
formulas using $H_{\nu}^{(1,2)}\equiv J_{\nu}\pm iY_{\nu}$, we obtain:
\begin{align}
\int_{0}^{1}ds\,H_{\nu}^{(1,2)}\left(\sigma s\right)s^{\nu} & =2^{\nu-1}\sigma^{-\nu}\sqrt{\pi}\,\Gamma\left(\nu+\frac{1}{2}\right)\left[H_{\nu}^{(1,2)}\left(\sigma\right)\mathbf{H}_{\nu-1}\left(\sigma\right)-H_{\nu-1}^{(1,2)}\left(\sigma\right)\mathbf{H}_{\nu}\left(\sigma\right)\right];\label{eq:intFormula1}\\
\int_{0}^{1}ds\,H_{\nu}^{(1,2)}\left(\sigma s\right)s^{\nu+1} & =\sigma^{-1}H_{\nu+1}^{(1,2)}\left(\sigma\right)\pm i2^{\nu+1}\sigma^{-\nu-2}\pi^{-1}\Gamma\left(\nu+1\right),\label{eq:intFormula2}
\end{align}
which give us the sought results for $\nu=0$: 
\begin{align}
\mathcal{I}_{0}\left(\sigma\right) & =\frac{\pi}{2}\sigma\left[H_{0}^{(1,2)}\left(\sigma\right)\mathbf{H}_{-1}\left(\sigma\right)+H_{1}^{(1,2)}\left(\sigma\right)\mathbf{H}_{0}\left(\sigma\right)\right];\label{eq:I0new}\\
\mathcal{I}_{1}\left(\sigma\right) & =\sigma H_{1}^{(1,2)}\left(\sigma\right)\pm\frac{2}{\pi}i.\label{eq:I1new}
\end{align}
It is fundamental to note that both $\mathcal{I}_{0}\left(\sigma\right)$
and $\mathcal{I}_{1}\left(\sigma\right)$ are regular functions with
removable singularity at $\sigma=0$, as displayed in Figure \ref{fig:profiles}.
 That is to say, the remaining integrals in (\ref{eq:I11new})-(\ref{eq:I22new})
are no longer singular and can be approximated with good accuracy
by standard Gauss-Legendre quadrature formulas. 

Let us now try to further simplify the resulting expressions.
\begin{align}
\int_{0}^{l_{n}}dx\,I_{\alpha}^{n}\left(x\right) & =\frac{1}{k}\int_{0}^{l_{n}}dx\,\mathcal{I}_{0}\left(kx\right)+\frac{1}{k}\int_{0}^{l_{n}}dx\,\mathcal{I}_{0}\left[k\left(l_{n}-x\right)\right]\nonumber \\
 & =\frac{1}{k^{2}}\int_{0}^{kl_{n}}dx\,\mathcal{I}_{0}\left(x\right)-\frac{1}{k^{2}}\int_{kl_{n}}^{0}dx\,\mathcal{I}_{0}\left(x\right)=\frac{2}{k^{2}}\int_{0}^{kl_{n}}dx\,\mathcal{I}_{0}\left(x\right);\label{eq:IalphaInt}
\end{align}
\begin{align}
\int_{0}^{l_{n}}dx\,I_{\beta}^{n}\left(x\right) & =\frac{1}{kl_{n}}\int_{0}^{l_{n}}dx\,x\mathcal{I}_{0}\left(kx\right)+\frac{1}{kl_{n}}\int_{0}^{l_{n}}dx\,x\mathcal{I}_{0}\left[k\left(l_{n}-x\right)\right]\nonumber \\
 & =\frac{1}{k^{3}l_{n}}\int_{0}^{kl_{n}}dx\,x\mathcal{I}_{0}\left(x\right)+\frac{1}{k^{2}l_{n}}\int_{0}^{kl_{n}}dx\,\left(l_{n}-\frac{x}{k}\right)\mathcal{I}_{0}\left(x\right)=\frac{1}{k^{2}}\int_{0}^{kl_{n}}dx\,\mathcal{I}_{0}\left(x\right).\label{eq:IbetaInt}
\end{align}
For the next integral, it is useful to note that:
\begin{align}
\int_{0}^{l_{n}}dx\,\left\{ \mathcal{I}_{1}\left[k\left(l_{n}-x\right)\right]-\mathcal{I}_{1}\left(kx\right)\right\}  & =\frac{1}{k}\int_{0}^{kl_{n}}dx\,\mathcal{I}_{1}\left(x\right)-\frac{1}{k}\int_{0}^{kl_{n}}dx\,\mathcal{I}_{1}\left(x\right)=0.
\end{align}
Therefore:
\begin{align}
\int_{0}^{l_{n}}dx\,I_{\gamma}^{n}\left(x\right) & =\int_{0}^{l_{n}}dx\,I_{\beta}^{n}\left(x\right)=\frac{1}{k^{2}}\int_{0}^{kl_{n}}dx\,\mathcal{I}_{0}\left(x\right).\label{eq:IgammaInt}
\end{align}
The last integral can be arrived at by considering the followings:
\begin{align}
\int_{0}^{l_{n}}dx\,\frac{x}{l_{n}}I_{\beta}^{n}\left(x\right) & =\frac{1}{kl_{n}^{2}}\int_{0}^{l_{n}}dx\,x^{2}\mathcal{I}_{0}\left(kx\right)+\frac{1}{kl_{n}^{2}}\int_{0}^{l_{n}}dx\,x^{2}\mathcal{I}_{0}\left[k\left(l_{n}-x\right)\right]\nonumber \\
 & =\frac{1}{k^{2}l_{n}^{2}}\int_{0}^{kl_{n}}dx\,\left[\frac{x^{2}}{k^{2}}\mathcal{I}_{0}\left(x\right)+\left(l_{n}-\frac{x}{k}\right)^{2}\mathcal{I}_{0}\left(x\right)\right]\nonumber \\
 & =\frac{1}{k^{2}l_{n}^{2}}\int_{0}^{kl_{n}}dx\,\left(l_{n}^{2}-\frac{2l_{n}}{k}x+\frac{2}{k^{2}}x^{2}\right)\mathcal{I}_{0}\left(x\right);
\end{align}
\begin{align}
 & \frac{1}{k^{2}l_{n}^{2}}\int_{0}^{l_{n}}dx\,\left\{ x\mathcal{I}_{1}\left[k\left(l_{n}-x\right)\right]-x\mathcal{I}_{1}\left(kx\right)\right\} \nonumber \\
 & =\frac{1}{k^{3}l_{n}^{2}}\int_{0}^{kl_{n}}dx\,\left[\left(l_{n}-\frac{x}{k}\right)\mathcal{I}_{1}\left(x\right)-\frac{x}{k}\mathcal{I}_{1}\left(x\right)\right]\nonumber \\
 & =\frac{1}{k^{3}l_{n}^{2}}\int_{0}^{kl_{n}}dx\,\left(l_{n}-\frac{2}{k}x\right)\mathcal{I}_{1}\left(x\right).
\end{align}
As the careful reader may notice, integrations involving $\mathcal{I}_{1}\left(x\right)$
and $x\,\mathcal{I}_{1}\left(x\right)$ can still be evaluated analytically
with the help of (\ref{eq:intFormula1}) and (\ref{eq:intFormula2}):
\begin{align}
\int_{0}^{kl_{n}}dx\,\mathcal{I}_{1}\left(x\right) & =\int_{0}^{kl_{n}}dx\,\left[xH_{1}^{(1,2)}\left(x\right)\pm\frac{2}{\pi}i\right]\nonumber \\
 & =\frac{kl_{n}\pi}{2}\left[H_{1}^{(1,2)}\left(kl_{n}\right)\mathbf{H}_{0}\left(kl_{n}\right)-H_{0}^{(1,2)}\left(kl_{n}\right)\mathbf{H}_{1}\left(kl_{n}\right)\right]\pm\frac{2kl_{n}}{\pi}i;
\end{align}
\begin{equation}
\int_{0}^{kl_{n}}dx\,x\mathcal{I}_{1}\left(x\right)=\int_{0}^{kl_{n}}dx\,\left[x^{2}H_{1}^{(1,2)}\left(x\right)\pm\frac{2i}{\pi}x\right]=k^{2}l_{n}^{2}H_{2}^{(1,2)}\left(kl_{n}\right)\pm\frac{i}{\pi}\left(k^{2}l_{n}^{2}+4\right),
\end{equation}
so that:
\begin{align}
\int_{0}^{l_{n}}dx\,I_{\delta}^{n}\left(x\right) & =\frac{1}{k^{2}l_{n}^{2}}\int_{0}^{kl_{n}}dx\,\left(l_{n}^{2}-\frac{2l_{n}}{k}x+\frac{2}{k^{2}}x^{2}\right)\mathcal{I}_{0}\left(x\right)+\nonumber \\
 & +\frac{\pi}{2k^{2}}\left[H_{1}^{(1,2)}\left(kl_{n}\right)\mathbf{H}_{0}\left(kl_{n}\right)-H_{0}^{(1,2)}\left(kl_{n}\right)\mathbf{H}_{1}\left(kl_{n}\right)\right]-\frac{2}{k^{2}}H_{2}^{(1,2)}\left(kl_{n}\right)\mp\frac{8i}{\pi k^{4}l_{n}^{2}}.\label{eq:IdeltaInt}
\end{align}
Gathering together the results in (\ref{eq:IalphaInt}), (\ref{eq:IbetaInt}),
(\ref{eq:IgammaInt}) and (\ref{eq:IdeltaInt}), we are finally able
to rewrite (\ref{eq:I11new})-(\ref{eq:I22new}) as follows:
\begin{align}
I_{11}^{n}=I_{22}^{n} & =\frac{i\pi}{8k^{2}}\left[H_{1}^{(1,2)}\left(kl_{n}\right)\mathbf{H}_{0}\left(kl_{n}\right)-H_{0}^{(1,2)}\left(kl_{n}\right)\mathbf{H}_{1}\left(kl_{n}\right)\right]+\nonumber \\
 & -\frac{i}{2k^{2}}H_{2}^{(1,2)}\left(kl_{n}\right)\pm\frac{2}{\pi k^{4}l_{n}^{2}}+\frac{i}{4k^{2}}\varGamma_{0}\left(kl_{n}\right)+\frac{i}{2k^{4}l_{n}^{2}}\varGamma_{2}\left(kl_{n}\right);\label{eq:I11I22}\\
I_{12}^{n}=I_{21}^{n} & =-I_{11}^{n}+\frac{i}{4k^{2}}\varGamma_{0}\left(kl_{n}\right),\label{eq:I12I21}
\end{align}
where the remaining integrals:
\begin{align}
\varGamma_{0}\left(\sigma\right) & \equiv\int_{0}^{\sigma}dx\,\mathcal{I}_{0}\left(x\right)=\frac{\pi}{2}\int_{0}^{\sigma}dx\,x\left[H_{0}^{(1,2)}\left(x\right)\mathbf{H}_{-1}\left(x\right)+H_{1}^{(1,2)}\left(x\right)\mathbf{H}_{0}\left(x\right)\right];\label{eq:gamma0}\\
\varGamma_{2}\left(\sigma\right) & \equiv\int_{0}^{\sigma}dx\,\mathcal{I}_{0}\left(x\right)\,x\left(x-\sigma\right)=\frac{\pi}{2}\int_{0}^{\sigma}dx\,x^{2}\left(x-\sigma\right)\left[H_{0}^{(1,2)}\left(x\right)\mathbf{H}_{-1}\left(x\right)+H_{1}^{(1,2)}\left(x\right)\mathbf{H}_{0}\left(x\right)\right]\label{eq:gamma2}
\end{align}
are left to Gauss-Legendre quadrature formulas.%

To summarize, the singular double integrals (\ref{eq:I11})-(\ref{eq:I22})
have been rewritten through (\ref{eq:I11I22}) and (\ref{eq:I12I21})
as a combination of the regular single integrals (\ref{eq:gamma0})-(\ref{eq:gamma2})
and of some well-known analytic functions.

\subsection{Integration over adjacent elements\label{subsec:adjacent}}

Whenever the elements $S_{m}$ and $S_{n}$ in (\ref{eq:singleLayer_Matrix})
are adjacent, the Green function (\ref{eq:Green}) diverges in correspondence
of the common vertex. As it is clear from Figure \ref{fig:combinations},
this can happen when the $i$-th and $j$-th mesh nodes coincide as
well as when they are first or second neighbors. In order to establish
the strength of the singularity, the following expansion of the zeroth
order Hankel function for small argument must be considered \cite{Abramowitz1972}:
\begin{equation}
H_{0}^{(1,2)}\left(z\right)\sim1\pm\frac{2i}{\pi}\left(\log z+\gamma-\log2\right),\label{eq:hankel0expansion}
\end{equation}
where $\gamma$ is the Euler's constant. It is important to note that
the integrand in (\ref{eq:singleLayer_Matrix}) is well-behaved when
the product of the functions $p_{i}^{m}\left(\mathbf{r}\right)$
and $p_{j}^{n}\left(\mathbf{r}'\right)$ is zero at the common
vertex, since it cancels the singularity of the Green function. Standard
Gauss-Legendre quadrature applies even when both $p_{i}^{m}\left(\mathbf{r}\right)$
and $p_{j}^{n}\left(\mathbf{r}'\right)$ are one at the common
vertex, i.e. for $i\equiv j$, as the logarithmic divergence is weak
enough to be integrable and the singular end point is not considered.
To improve the accuracy of the results, adaptive integration algorithms
can be used in this case (see, for instance, \cite{Ramdas2002}).

\section{Double layer and adjoint double layer operators\label{sec:doubleLayer}}

\subsection{Integration over coincident elements}

When $S_{m}=S_{n}$, the double integrals appearing in (\ref{eq:doubleLayer_Matrix})
and (\ref{eq:adjointDoubleLayer_Matrix}) vanish identically, as follows
from the fact that $\mathbf{R}\perp\mathbf{n}$ when both
$\mathbf{r}$ and $\mathbf{r}'$ lie on the same segment with
unit normal $\mathbf{n}$.

\subsection{Integration over adjacent elements\label{subsec:adjacentDouble}}

Making use of the following expansion of the first order Hankel function
for small argument \cite{Abramowitz1972}:
\begin{equation}
H_{1}^{(1,2)}\left(z\right)\sim\mp\frac{2i}{\pi z},\label{eq:HankelExpansion}
\end{equation}
it is easy to see that the only singular contribution to the integrations
over adjacent segments in (\ref{eq:doubleLayer_Matrix}) and (\ref{eq:adjointDoubleLayer_Matrix})
is achieved, once again, when both $p_{i}^{m}\left(\mathbf{r}\right)$
and $p_{j}^{n}\left(\mathbf{r}'\right)$ are one at the common
vertex. In this case, however, the singularity is stronger than that
in Section \ref{subsec:adjacent} and direct use of Gauss-Legendre
quadrature formulas may lead to inaccurate results. In order to avoid
this issue, a viable technique consists in introducing a coordinate
transformation with vanishing Jacobian at the common vertex to cancel
the singularity \cite{Sutradhar2008}. In \ref{sec:appendix},
this method is applied to the double layer integral; of course, the
same procedure can be used in the adjoint double layer case with the
appropriate changes.

\section{Hypersingular operator: direct method\label{sec:directMethod}}

\subsection{Integration over coincident elements\label{subsec:coincidentHypersingular}}

When $S_{m}=S_{n}$, the singular double integrals in (\ref{eq:hyperSingular_Matrix})
are of the form:
\begin{equation}
\frac{ik}{4}\int_{S_{n}}d\mathbf{r}\int_{S_{n}}d\mathbf{r}'\frac{H_{1}^{(1,2)}\left(k\left|\mathbf{r}-\mathbf{r}'\right|\right)}{\left|\mathbf{r}-\mathbf{r}'\right|}p_{i}^{n}\left(\mathbf{r}\right)p_{j}^{n}\left(\mathbf{r}'\right).
\end{equation}
Adopting the same convention introduced in Section \ref{subsec:coincidentIntSingleLayer},
we have:
\begin{equation}
\frac{ik}{4}\int_{0}^{l_{n}}dx\int_{-x}^{l_{n}-x}dx'\frac{H_{1}^{(1,2)}\left(k\left|x'\right|\right)}{\left|x'\right|}p_{i}^{n}\left(x\right)\widetilde{p}_{j}^{n}\left(x',x\right).
\end{equation}
Making use of the recurrence relations for Hankel functions \cite{Gradshteyn2007}:
\begin{align}
zH_{\nu-1}^{(1,2)}\left(z\right)+zH_{\nu+1}^{(1,2)}\left(z\right) & =2\nu H_{\nu}^{(1,2)}\left(z\right);\\
H_{\nu-1}^{(1,2)}\left(z\right)-H_{\nu+1}^{(1,2)}\left(z\right) & =2\frac{d}{dz}H_{\nu}^{(1,2)}\left(z\right)\equiv2\dot{H}_{\nu}^{(1,2)}\left(z\right),\label{eq:recurrenceDerivative}
\end{align}
we get: 
\begin{equation}
\frac{ik^{2}}{4}\int_{0}^{l_{n}}dx\int_{-x}^{l_{n}-x}dx'\left[H_{0}^{(1,2)}\left(k\left|x'\right|\right)-\dot{H}_{1}^{(1,2)}\left(k\left|x'\right|\right)\right]p_{i}^{n}\left(x\right)\widetilde{p}_{j}^{n}\left(x',x\right).
\end{equation}
This last equation, together with (\ref{eq:bfi}), (\ref{eq:bfj})
and (\ref{eq:I11})-(\ref{eq:I22}), gives rise to four expressions:
\begin{align}
\varUpsilon_{11}^{n} & \equiv k^{2}I_{11}^{n}-\frac{ik^{2}}{4}\int_{0}^{l_{n}}dx\int_{-x}^{l_{n}-x}dx'\dot{H}_{1}^{(1,2)}\left(k\left|x'\right|\right)\left(1-\frac{x}{l_{n}}\right)\left[1-\frac{\left(x'+x\right)}{l_{n}}\right];\label{eq:U11}\\
\varUpsilon_{12}^{n} & \equiv k^{2}I_{12}^{n}-\frac{ik^{2}}{4}\int_{0}^{l_{n}}dx\int_{-x}^{l_{n}-x}dx'\dot{H}_{1}^{(1,2)}\left(k\left|x'\right|\right)\left(1-\frac{x}{l_{n}}\right)\left[\frac{\left(x'+x\right)}{l_{n}}\right];\\
\varUpsilon_{21}^{n} & \equiv k^{2}I_{21}^{n}-\frac{ik^{2}}{4}\int_{0}^{l_{n}}dx\int_{-x}^{l_{n}-x}dx'\dot{H}_{1}^{(1,2)}\left(k\left|x'\right|\right)\left(\frac{x}{l_{n}}\right)\left[1-\frac{\left(x'+x\right)}{l_{n}}\right];\\
\varUpsilon_{22}^{n} & \equiv k^{2}I_{22}^{n}-\frac{ik^{2}}{4}\int_{0}^{l_{n}}dx\int_{-x}^{l_{n}-x}dx'\dot{H}_{1}^{(1,2)}\left(k\left|x'\right|\right)\left(\frac{x}{l_{n}}\right)\left[\frac{\left(x'+x\right)}{l_{n}}\right],\label{eq:U22}
\end{align}
which can be recast in the following form:
\begin{align}
\varUpsilon_{11}^{n} & =k^{2}I_{11}^{n}-\frac{ik^{2}}{4}\int_{0}^{l_{n}}dx\left[\varUpsilon_{\alpha}^{n}\left(x\right)-\varUpsilon_{\beta}^{n}\left(x\right)-\varUpsilon_{\gamma}^{n}\left(x\right)+\varUpsilon_{\delta}^{n}\left(x\right)\right];\label{eq:U11new}\\
\varUpsilon_{12}^{n} & =k^{2}I_{12}^{n}-\frac{ik^{2}}{4}\int_{0}^{l_{n}}dx\left[\varUpsilon_{\gamma}^{n}\left(x\right)-\varUpsilon_{\delta}^{n}\left(x\right)\right];\\
\varUpsilon_{21}^{n} & =k^{2}I_{12}^{n}-\frac{ik^{2}}{4}\int_{0}^{l_{n}}dx\left[\varUpsilon_{\beta}^{n}\left(x\right)-\varUpsilon_{\delta}^{n}\left(x\right)\right];\\
\varUpsilon_{22}^{n} & =k^{2}I_{11}^{n}-\frac{ik^{2}}{4}\int_{0}^{l_{n}}dx\,\varUpsilon_{\delta}^{n}\left(x\right),\label{eq:U22new}
\end{align}
where (\ref{eq:I11I22}), (\ref{eq:I12I21}) have been used and:
\begin{align}
\varUpsilon_{\alpha}^{n}\left(x\right) & \equiv\int_{-x}^{l_{n}-x}dx'\dot{H}_{1}^{(1,2)}\left(k\left|x'\right|\right);\\
\varUpsilon_{\beta}^{n}\left(x\right) & \equiv\frac{x}{l_{n}}\int_{-x}^{l_{n}-x}dx'\dot{H}_{1}^{(1,2)}\left(k\left|x'\right|\right);\\
\varUpsilon_{\gamma}^{n}\left(x\right) & \equiv\int_{-x}^{l_{n}-x}dx'\dot{H}_{1}^{(1,2)}\left(k\left|x'\right|\right)\frac{\left(x'+x\right)}{l_{n}};\\
\varUpsilon_{\delta}^{n}\left(x\right) & \equiv\frac{x}{l_{n}}\int_{-x}^{l_{n}-x}dx'\dot{H}_{1}^{(1,2)}\left(k\left|x'\right|\right)\frac{\left(x'+x\right)}{l_{n}}.
\end{align}
As in Section \ref{subsec:coincidentIntSingleLayer}, these integrals
will be split in order to get rid of absolute values. For instance:
\begin{align}
\varUpsilon_{\alpha}^{n}\left(x\right) & =\int_{-x}^{0}dx'\dot{H}_{1}^{(1,2)}\left(-kx'\right)+\int_{0}^{l_{n}-x}dx'\dot{H}_{1}^{(1,2)}\left(kx'\right)\nonumber \\
 & =\frac{1}{k}\int_{0}^{kx}dx'\dot{H}_{1}^{(1,2)}\left(x'\right)+\frac{1}{k}\int_{0}^{k\left(l_{n}-x\right)}dx'\dot{H}_{1}^{(1,2)}\left(x'\right)
\end{align}
and similarly for the other three integrals. Then, we define:
\begin{align}
\varUpsilon_{0}\left(\sigma\right) & \equiv\int_{0}^{\sigma}dx'\dot{H}_{1}^{(1,2)}\left(x'\right);\label{eq:U0}\\
\varUpsilon_{1}\left(\sigma\right) & \equiv\int_{0}^{\sigma}dx'\dot{H}_{1}^{(1,2)}\left(x'\right)\,x',\label{eq:U1}
\end{align}
so that:
\begin{align}
\varUpsilon_{\alpha}^{n}\left(x\right) & =\frac{1}{k}\left\{ \varUpsilon_{0}\left(kx\right)+\varUpsilon_{0}\left[k\left(l_{n}-x\right)\right]\right\} ;\label{eq:Ualpha}\\
\varUpsilon_{\beta}^{n}\left(x\right) & =\frac{x}{l_{n}}\varUpsilon_{\alpha}^{n}\left(x\right);\\
\varUpsilon_{\gamma}^{n}\left(x\right) & =\varUpsilon_{\beta}^{n}\left(x\right)+\frac{1}{k^{2}l_{n}}\left\{ \varUpsilon_{1}\left[k\left(l_{n}-x\right)\right]-\varUpsilon_{1}\left(kx\right)\right\} ;\\
\varUpsilon_{\delta}^{n}\left(x\right) & =\frac{x}{l_{n}}\varUpsilon_{\gamma}^{n}\left(x\right).\label{eq:Udelta}
\end{align}
Unfortunately, both (\ref{eq:U0}) and (\ref{eq:U1}) are singular.
A regularization for $\varUpsilon_{0}\left(\sigma\right)$ can be
achieved by analytic continuation:
\begin{align}
\varUpsilon_{0}\left(\sigma\right) & =\frac{1}{2}\lim_{\varepsilon\rightarrow0}\left[\int_{0}^{\sigma}dx'\dot{H}_{1}^{(1,2)}\left(x'-i\varepsilon\right)+\int_{0}^{\sigma}dx'\dot{H}_{1}^{(1,2)}\left(x'+i\varepsilon\right)\right]\nonumber \\
 & =\frac{1}{2}\lim_{\varepsilon\rightarrow0}\left[\int_{-i\varepsilon}^{\sigma-i\varepsilon}dx'\dot{H}_{1}^{(1,2)}\left(x'\right)+\int_{i\varepsilon}^{\sigma+i\varepsilon}dx'\dot{H}_{1}^{(1,2)}\left(x'\right)\right]\nonumber \\
 & =\frac{1}{2}\lim_{\varepsilon\rightarrow0}\left[H_{1}^{(1,2)}\left(\sigma-i\varepsilon\right)-H_{1}^{(1,2)}\left(-i\varepsilon\right)+H_{1}^{(1,2)}\left(\sigma+i\varepsilon\right)-H_{1}^{(1,2)}\left(i\varepsilon\right)\right]\nonumber \\
 & =H_{1}^{(1,2)}\left(\sigma\right)-\frac{1}{2}\lim_{\varepsilon\rightarrow0}\left[H_{1}^{(1,2)}\left(i\varepsilon\right)+H_{1}^{(1,2)}\left(-i\varepsilon\right)\right]=H_{1}^{(1,2)}\left(\sigma\right),
\end{align}
where the expansion (\ref{eq:HankelExpansion}) has been employed
in the last step. Now, exploiting the fact that the divergent integral
(\ref{eq:U1}) only appears through the difference $\varUpsilon_{1}\left[k\left(l_{n}-x\right)\right]-\varUpsilon_{1}\left(kx\right)$,
let:
\begin{align}
\delta\varUpsilon_{1}\left(\sigma_{a},\sigma_{b}\right) & \equiv\varUpsilon_{1}\left(\sigma_{b}\right)-\varUpsilon_{1}\left(\sigma_{a}\right)=\fint_{\sigma_{a}}^{\sigma_{b}}dx'\dot{H}_{1}^{(1,2)}\left(x'\right)\,x'\nonumber \\
 & =\lim_{\varepsilon\rightarrow0}\left[\int_{\varepsilon}^{\sigma_{b}}dx'\dot{H}_{1}^{(1,2)}\left(x'\right)\,x'-\int_{\varepsilon}^{\sigma_{a}}dx'\dot{H}_{1}^{(1,2)}\left(x'\right)\,x'\right].
\end{align}
Each of the two terms in the last expression can be integrated by
parts:
\begin{align}
\int_{\varepsilon}^{\sigma}dx'\dot{H}_{1}^{(1,2)}\left(x'\right)\,x' & =\sigma H_{1}^{(1,2)}\left(\sigma\right)-\varepsilon H_{1}^{(1,2)}\left(\varepsilon\right)-\int_{\varepsilon}^{\sigma}dx'H_{1}^{(1,2)}\left(x'\right)
\end{align}
where, using (\ref{eq:recurrenceDerivative}):
\begin{align}
\int_{\varepsilon}^{\sigma}dx'H_{1}^{(1,2)}\left(x'\right) & =-\int_{\varepsilon}^{\sigma}dx'\dot{H}_{0}^{(1,2)}\left(x'\right)=-\left[H_{0}^{(1,2)}\left(\sigma\right)-H_{0}^{(1,2)}\left(\varepsilon\right)\right].\label{eq:H1integral}
\end{align}
Then, taking the difference:
\begin{align}
\delta\varUpsilon_{1}\left(\sigma_{a},\sigma_{b}\right) & =\sigma_{b}H_{1}^{(1,2)}\left(\sigma_{b}\right)-\sigma_{a}H_{1}^{(1,2)}\left(\sigma_{a}\right)+H_{0}^{(1,2)}\left(\sigma_{b}\right)-H_{0}^{(1,2)}\left(\sigma_{a}\right).
\end{align}

Let us apply the above results to simplify expressions (\ref{eq:Ualpha})-(\ref{eq:Udelta}):
\begin{align}
\varUpsilon_{\alpha}^{n}\left(x\right) & =\frac{1}{k}\left\{ H_{1}^{(1,2)}\left(kx\right)+H_{1}^{(1,2)}\left[k\left(l_{n}-x\right)\right]\right\} ;\\
\varUpsilon_{\beta}^{n}\left(x\right) & =\frac{x}{l_{n}}\varUpsilon_{\alpha}^{n}\left(x\right);\\
\varUpsilon_{\gamma}^{n}\left(x\right) & =\varUpsilon_{\beta}^{n}\left(x\right)+\frac{1}{k^{2}l_{n}}\left\{ k\left(l_{n}-x\right)H_{1}^{(1,2)}\left[k\left(l_{n}-x\right)\right]-kxH_{1}^{(1,2)}\left(kx\right)+\right.\\
 & \left.+H_{0}^{(1,2)}\left[k\left(l_{n}-x\right)\right]-H_{0}^{(1,2)}\left(kx\right)\right\} ;\\
\varUpsilon_{\delta}^{n}\left(x\right) & =\frac{x}{l_{n}}\varUpsilon_{\gamma}^{n}\left(x\right).
\end{align}
In order to be able to rewrite (\ref{eq:U11new})-(\ref{eq:U22new})
in closed form, the previous equations will be integrated analytically
between $\varepsilon$ and $l_{n}-\varepsilon$ and the limit $\varepsilon\rightarrow0$
will be taken explicitly whenever possible, otherwise implicitly assumed.
Starting with $\varUpsilon_{\alpha}^{n}\left(x\right)$, making use
of the changes of variable $x\rightarrow kx$, $x\rightarrow k\left(l_{n}-x\right)$
and of formula (\ref{eq:H1integral}), we have:
\begin{align}
\int_{\varepsilon}^{l_{n}-\varepsilon}dx\,\varUpsilon_{\alpha}^{n}\left(x\right) & =\frac{1}{k}\int_{\varepsilon}^{l_{n}-\varepsilon}dx\,\left\{ H_{1}^{(1,2)}\left(kx\right)+H_{1}^{(1,2)}\left[k\left(l_{n}-x\right)\right]\right\} \nonumber \\
 & =\frac{1}{k^{2}}\int_{k\varepsilon}^{kl_{n}}dx\,H_{1}^{(1,2)}\left(x\right)-\frac{1}{k^{2}}\int_{kl_{n}}^{k\varepsilon}dx\,H_{1}^{(1,2)}\left(x\right)\nonumber \\
 & =\frac{2}{k^{2}}\int_{k\varepsilon}^{kl_{n}}dx\,H_{1}^{(1,2)}\left(x\right)=\frac{2}{k^{2}}\left[H_{0}^{(1,2)}\left(k\varepsilon\right)-H_{0}^{(1,2)}\left(kl_{n}\right)\right].\label{eq:UalphaInt}
\end{align}
The integral of $\varUpsilon_{\beta}^{n}\left(x\right)$ is treated
analogously: 
\begin{align}
\int_{\varepsilon}^{l_{n}-\varepsilon}dx\,\varUpsilon_{\beta}^{n}\left(x\right) & =\frac{1}{kl_{n}}\int_{\varepsilon}^{l_{n}-\varepsilon}dx\,\left\{ xH_{1}^{(1,2)}\left(kx\right)+xH_{1}^{(1,2)}\left[k\left(l_{n}-x\right)\right]\right\} \nonumber \\
 & =\frac{1}{k^{3}l_{n}}\int_{k\varepsilon}^{kl_{n}}dx\,xH_{1}^{(1,2)}\left(x\right)-\frac{1}{k^{2}l_{n}}\int_{kl_{n}}^{k\varepsilon}dx\,\left(l_{n}-\frac{x}{k}\right)H_{1}^{(1,2)}\left(x\right)\nonumber \\
 & =\frac{1}{k^{3}l_{n}}\int_{k\varepsilon}^{kl_{n}}dx\,xH_{1}^{(1,2)}\left(x\right)+\frac{1}{k^{2}}\int_{k\varepsilon}^{kl_{n}}dx\,H_{1}^{(1,2)}\left(x\right)-\frac{1}{k^{3}l_{n}}\int_{k\varepsilon}^{kl_{n}}dx\,xH_{1}^{(1,2)}\left(x\right)\nonumber \\
 & =\frac{1}{k^{2}}\int_{k\varepsilon}^{kl_{n}}dx\,H_{1}^{(1,2)}\left(x\right)=\frac{1}{k^{2}}\left[H_{0}^{(1,2)}\left(k\varepsilon\right)-H_{0}^{(1,2)}\left(kl_{n}\right)\right].\label{eq:UbetaInt}
\end{align}
To compute the integral of $\varUpsilon_{\gamma}^{n}\left(x\right)$,
we notice that: 
\begin{align}
 & \frac{1}{k^{2}l_{n}}\int_{\varepsilon}^{l_{n}-\varepsilon}dx\,\left\{ k\left(l_{n}-x\right)H_{1}^{(1,2)}\left[k\left(l_{n}-x\right)\right]-kxH_{1}^{(1,2)}\left(kx\right)\right\} \nonumber \\
 & =\frac{1}{k^{3}l_{n}}\left[\int_{k\varepsilon}^{kl_{n}}dx\,xH_{1}^{(1,2)}\left(x\right)-\int_{k\varepsilon}^{kl_{n}}dx\,xH_{1}^{(1,2)}\left(x\right)\right]=0
\end{align}
and similarly:
\begin{align}
 & \frac{1}{k^{2}l_{n}}\int_{\varepsilon}^{l_{n}-\varepsilon}dx\,\left\{ H_{0}^{(1,2)}\left[k\left(l_{n}-x\right)\right]-H_{0}^{(1,2)}\left(kx\right)\right\} \nonumber \\
 & =\frac{1}{k^{3}l_{n}}\left[\int_{k\varepsilon}^{kl_{n}}dx\,H_{0}^{(1,2)}\left(x\right)-\int_{k\varepsilon}^{kl_{n}}dx\,H_{0}^{(1,2)}\left(x\right)\right]=0.
\end{align}
Therefore:
\begin{equation}
\int_{\varepsilon}^{l_{n}-\varepsilon}dx\,\varUpsilon_{\gamma}^{n}\left(x\right)=\int_{\varepsilon}^{l_{n}-\varepsilon}dx\,\varUpsilon_{\beta}^{n}\left(x\right)=\frac{1}{k^{2}}\left[H_{0}^{(1,2)}\left(k\varepsilon\right)-H_{0}^{(1,2)}\left(kl_{n}\right)\right].\label{eq:UgammaInt}
\end{equation}
The following results prove useful for the evaluation of the last
integral:
\begin{align}
\int_{\varepsilon}^{l_{n}-\varepsilon}dx\,\frac{x}{l_{n}}\varUpsilon_{\beta}^{n}\left(x\right) & =\frac{1}{kl_{n}^{2}}\int_{\varepsilon}^{l_{n}-\varepsilon}dx\,\left\{ x^{2}H_{1}^{(1,2)}\left(kx\right)+x^{2}H_{1}^{(1,2)}\left[k\left(l_{n}-x\right)\right]\right\} \nonumber \\
 & =\frac{1}{k^{2}l_{n}^{2}}\int_{k\varepsilon}^{kl_{n}}dx\,\left[\frac{x^{2}}{k^{2}}H_{1}^{(1,2)}\left(x\right)+\left(l_{n}-\frac{x}{k}\right)^{2}H_{1}^{(1,2)}\left(x\right)\right]\nonumber \\
 & =\frac{1}{k^{2}l_{n}^{2}}\int_{k\varepsilon}^{kl_{n}}dx\,\left(l_{n}^{2}-\frac{2l_{n}}{k}x+\frac{2}{k^{2}}x^{2}\right)H_{1}^{(1,2)}\left(x\right);
\end{align}
\begin{align}
 & \frac{1}{k^{2}l_{n}^{2}}\int_{\varepsilon}^{l_{n}-\varepsilon}dx\,\left\{ k\left(l_{n}-x\right)xH_{1}^{(1,2)}\left[k\left(l_{n}-x\right)\right]-kx^{2}H_{1}^{(1,2)}\left(kx\right)\right\} \nonumber \\
 & =\frac{1}{k^{3}l_{n}^{2}}\int_{k\varepsilon}^{kl_{n}}dx\,\left[x\left(l_{n}-\frac{x}{k}\right)H_{1}^{(1,2)}\left(x\right)-\frac{x^{2}}{k}H_{1}^{(1,2)}\left(x\right)\right]\nonumber \\
 & =\frac{1}{k^{2}l_{n}^{2}}\int_{k\varepsilon}^{kl_{n}}dx\,\left(\frac{l_{n}}{k}x-\frac{2}{k^{2}}x^{2}\right)H_{1}^{(1,2)}\left(x\right);
\end{align}
\begin{align}
 & \frac{1}{k^{2}l_{n}^{2}}\int_{\varepsilon}^{l_{n}-\varepsilon}dx\,\left\{ xH_{0}^{(1,2)}\left[k\left(l_{n}-x\right)\right]-xH_{0}^{(1,2)}\left(kx\right)\right\} \nonumber \\
 & =\frac{1}{k^{3}l_{n}^{2}}\int_{k\varepsilon}^{kl_{n}}dx\,\left[\left(l_{n}-\frac{x}{k}\right)H_{0}^{(1,2)}\left(x\right)-\frac{x}{k}H_{0}^{(1,2)}\left(x\right)\right]\nonumber \\
 & =\frac{1}{k^{3}l_{n}}\int_{k\varepsilon}^{kl_{n}}dx\,H_{0}^{(1,2)}\left(x\right)-\frac{2}{k^{4}l_{n}^{2}}\int_{k\varepsilon}^{kl_{n}}dx\,xH_{0}^{(1,2)}\left(x\right)\nonumber \\
 & =\frac{1}{k^{3}l_{n}}\mathcal{I}_{0}\left(kl_{n}\right)-\frac{2}{k^{4}l_{n}^{2}}\mathcal{I}_{1}\left(kl_{n}\right),
\end{align}
where (\ref{eq:I0}) and (\ref{eq:I1}) have been employed in the
last step. Combining the previous expressions, we obtain:
\begin{align}
\int_{\varepsilon}^{l_{n}-\varepsilon}dx\,\varUpsilon_{\delta}^{n}\left(x\right) & =\frac{1}{k^{2}l_{n}^{2}}\int_{k\varepsilon}^{kl_{n}}dx\,\left(l_{n}^{2}-\frac{l_{n}}{k}x\right)H_{1}^{(1,2)}\left(x\right)+\frac{1}{k^{3}l_{n}}\mathcal{I}_{0}\left(kl_{n}\right)-\frac{2}{k^{4}l_{n}^{2}}\mathcal{I}_{1}\left(kl_{n}\right)\nonumber \\
 & =\frac{1}{k^{2}}\int_{k\varepsilon}^{kl_{n}}dx\,H_{1}^{(1,2)}\left(x\right)-\frac{1}{k^{3}l_{n}}\int_{k\varepsilon}^{kl_{n}}dx\,xH_{1}^{(1,2)}\left(x\right)+\frac{1}{k^{3}l_{n}}\mathcal{I}_{0}\left(kl_{n}\right)-\frac{2}{k^{4}l_{n}^{2}}\mathcal{I}_{1}\left(kl_{n}\right),
\end{align}
which can be simplified using (\ref{eq:H1integral}), the following
formula:
\begin{align}
\int_{0}^{\sigma}dx\,H_{1}^{(1,2)}\left(x\right)\,x & =\frac{\pi\sigma}{2}\left[H_{1}^{(1,2)}\left(\sigma\right)\mathbf{H}_{0}\left(\sigma\right)-H_{0}^{(1,2)}\left(\sigma\right)\mathbf{H}_{1}\left(\sigma\right)\right],
\end{align}
derived from (\ref{eq:intFormula1}) with $\nu=1$, and equations
(\ref{eq:I0new}) and (\ref{eq:I1new}):
\begin{align}
\int_{\varepsilon}^{l_{n}-\varepsilon}dx\,\varUpsilon_{\delta}^{n}\left(x\right) & =\frac{1}{k^{2}}\left[H_{0}^{(1,2)}\left(k\varepsilon\right)-H_{0}^{(1,2)}\left(kl_{n}\right)\right]+\nonumber \\
 & +\frac{\pi}{2k^{2}}\left[H_{0}^{(1,2)}\left(kl_{n}\right)\mathbf{H}_{1}\left(kl_{n}\right)-H_{1}^{(1,2)}\left(kl_{n}\right)\mathbf{H}_{0}\left(kl_{n}\right)\right]+\nonumber \\
 & +\frac{\pi}{2k^{2}}\left[H_{0}^{(1,2)}\left(kl_{n}\right)\mathbf{H}_{-1}\left(kl_{n}\right)+H_{1}^{(1,2)}\left(kl_{n}\right)\mathbf{H}_{0}\left(kl_{n}\right)\right]+\nonumber \\
 & -\frac{2}{k^{4}l_{n}^{2}}\left[kl_{n}H_{1}^{(1,2)}\left(kl_{n}\right)\pm\frac{2}{\pi}i\right].
\end{align}
Finally, by the recurrence relations for Struve functions \cite{Gradshteyn2007}:
\begin{align}
\mathbf{H}_{\nu-1}\left(z\right) & +\mathbf{H}_{\nu+1}\left(z\right)=\frac{2\nu}{z}\mathbf{H}_{\nu}\left(z\right)+\frac{1}{\sqrt{\pi}}\left(\frac{z}{2}\right)^{\nu}\frac{1}{\Gamma\left(\nu+\frac{3}{2}\right)},
\end{align}
we get:
\begin{align}
\int_{\varepsilon}^{l_{n}-\varepsilon}dx\,\varUpsilon_{\delta}^{n}\left(x\right) & =\frac{1}{k^{2}}H_{0}^{(1,2)}\left(k\varepsilon\right)-\frac{2}{k^{4}l_{n}^{2}}\left[kl_{n}H_{1}^{(1,2)}\left(kl_{n}\right)\pm\frac{2}{\pi}i\right].\label{eq:UdeltaInt}
\end{align}

Replacing (\ref{eq:UalphaInt}), (\ref{eq:UbetaInt}), (\ref{eq:UgammaInt})
and (\ref{eq:UdeltaInt}) into equations (\ref{eq:U11new})-(\ref{eq:U22new})
leads to:
\begin{align}
\varUpsilon_{11}^{n} & =\varUpsilon_{22}^{n}=k^{2}I_{11}^{n}-\frac{i}{4}H_{0}^{(1,2)}\left(k\varepsilon\right)+\frac{i}{2kl_{n}}H_{1}^{(1,2)}\left(kl_{n}\right)\mp\frac{1}{\pi k^{2}l_{n}^{2}};\label{eq:upsilon11and22}\\
\varUpsilon_{12}^{n} & =\varUpsilon_{21}^{n}=k^{2}I_{12}^{n}+\frac{i}{4}H_{0}^{(1,2)}\left(kl_{n}\right)-\frac{i}{2kl_{n}}H_{1}^{(1,2)}\left(kl_{n}\right)\pm\frac{1}{\pi k^{2}l_{n}^{2}},\label{eq:upsilon12and21}
\end{align}
where $I_{11}^{n}$ and $I_{12}^{n}$ are given by (\ref{eq:I11I22})
and (\ref{eq:I12I21}), respectively, and the limit $\varepsilon\rightarrow0$
is assumed. Owing to the logarithmic divergence of the function $H_{0}^{(1,2)}\left(k\varepsilon\right)$,
both $\varUpsilon_{11}^{n}$ and $\varUpsilon_{22}^{n}$ are singular.
In order for the matrix elements $N_{ii}$ in (\ref{eq:hyperSingular_Matrix})
to be well-defined, these singularities as well as those arising from
the integration over adjacent segments must cancel in the sum over
$\left\{ m\in i\right\} \times\left\{ n\in i\right\} $ (left column
of Figure \ref{fig:combinations}). This fundamental condition is
checked in the next subsection. By explicitly removing the divergent
term, equation (\ref{eq:upsilon11and22}) becomes:
\begin{align}
\varUpsilon_{11}^{n}\mp\frac{1}{2\pi}\log\varepsilon & =k^{2}I_{11}^{n}+\frac{i}{2kl_{n}}H_{1}^{(1,2)}\left(kl_{n}\right)-\frac{i}{4}\pm\frac{1}{2\pi}\left(\gamma+\log\frac{k}{2}-\frac{2}{k^{2}l_{n}^{2}}\right),\label{eq:ups11and22regularized}
\end{align}
where (\ref{eq:hankel0expansion}) has been used.%

\subsection{Integration over adjacent elements\label{subsec:adjacentHypersingular}}

Let us consider the following expansions for small $z$ \cite{Abramowitz1972}:
\begin{align}
\frac{H_{1}^{(1,2)}\left(z\right)}{z} & \sim\frac{1}{2}\mp\frac{2i}{\pi}\frac{1}{z^{2}}\pm\frac{i}{\pi}\left(\log z-\frac{1}{2}+\gamma-\log2\right);\label{eq:H1overZ}\\
H_{2}^{(1,2)}\left(z\right) & \sim\mp\frac{i}{\pi}\mp\frac{4i}{\pi}\frac{1}{z^{2}}.\label{eq:H2}
\end{align}
When only one of the two functions $p_{i}^{m}\left(\mathbf{r}\right)$
and $p_{j}^{n}\left(\mathbf{r}'\right)$ equals unity at the common
vertex, the integrations over adjacent segments in (\ref{eq:hyperSingular_Matrix})
can be treated via singularity cancellation as it is done in
\ref{sec:appendix} for the double layer case. Whereas the logarithmic
term does not constitute an issue, as explained in Section \ref{subsec:adjacent},
it is apparent that a non-integrable singularity of the form $\sim R^{-2}$
arises when both basis functions are one at the common vertex. Once
properly isolated, such divergent term proves to be equal and opposite
to that appearing in (\ref{eq:upsilon11and22}), as it is shown below.

With reference to (\ref{eq:H1overZ}) and (\ref{eq:H2}), a regularized
version of the first and second order Hankel functions can be defined
by singularity subtraction:
\begin{align}
\widetilde{H}_{1}^{(1,2)}\left(z\right) & \equiv H_{1}^{(1,2)}\left(z\right)\pm\frac{2i}{\pi}\frac{1}{z};\label{eq:regularH1}\\
\widetilde{H}_{2}^{(1,2)}\left(z\right) & \equiv H_{2}^{(1,2)}\left(z\right)\pm\frac{4i}{\pi}\frac{1}{z^{2}}.\label{eq:regularH2}
\end{align}
This allows us to rewrite the divergent adjacent integrations in (\ref{eq:hyperSingular_Matrix})
as the sum of a regular part and a singular part:
\begin{equation}
\int_{S_{m}}d\mathbf{r}\int_{S_{n}}d\mathbf{r}'\frac{\partial^{2}g\left(\mathbf{r},\mathbf{r}'\right)}{\partial n\partial n'}p_{i}^{m}\left(\mathbf{r}\right)p_{j}^{n}\left(\mathbf{r}'\right)=\varUpsilon_{reg}^{m,n}+\varUpsilon_{sing}^{m,n},
\end{equation}
where:
\begin{align}
\varUpsilon_{reg}^{m,n} & =\frac{ik^{2}}{4}\int_{S_{m}}d\mathbf{r}\int_{S_{n}}d\mathbf{r}'\left[\frac{\widetilde{H}_{1}^{(1,2)}\left(kR\right)}{kR}\left(\mathbf{n}\cdot\mathbf{n}'\right)-\widetilde{H}_{2}^{(1,2)}\left(kR\right)\frac{\left(\mathbf{R}\cdot\mathbf{n}\right)\left(\mathbf{R}\cdot\mathbf{n}'\right)}{R^{2}}\right]p_{i}^{m}\left(\mathbf{r}\right)p_{j}^{n}\left(\mathbf{r}'\right);\label{eq:regInt}\\
\varUpsilon_{sing}^{m,n} & =\pm\frac{1}{2\pi}\int_{S_{m}}d\mathbf{r}\int_{S_{n}}d\mathbf{r}'\left[\frac{\left(\mathbf{n}\cdot\mathbf{n}'\right)}{R^{2}}-2\frac{\left(\mathbf{R}\cdot\mathbf{n}\right)\left(\mathbf{R}\cdot\mathbf{n}'\right)}{R^{4}}\right]p_{i}^{m}\left(\mathbf{r}\right)p_{j}^{n}\left(\mathbf{r}'\right)\label{eq:singInt}
\end{align}
and basis functions that equal one at the common vertex are assumed.
The regular integral $\varUpsilon_{reg}^{m,n}$ is left to Gauss-Legendre
quadrature formulas. Conversely, by fixing the reference frames $O'x'y'$,
$Oxy$ as in Figure \ref{fig:referenceFrames} (right side), so that
the basis functions are given by:
\begin{equation}
p_{i}^{m}\left(x\right)\equiv\begin{cases}
1-\frac{x}{l_{m}} & \mathrm{if}\:x_{i}=0;\\
\frac{x}{l_{m}} & \mathrm{if}\:x_{i}=l_{m},
\end{cases}\label{eq:pim}
\end{equation}
\begin{equation}
p_{j}^{n}\left(x',x\right)\equiv\begin{cases}
1-\frac{\left(x'+x\cos\theta\right)}{l_{n}} & \mathrm{if}\:x_{j}'=-x\cos\theta;\\
\frac{\left(x'+x\cos\theta\right)}{l_{n}} & \mathrm{if}\:x_{j}'=l_{n}-x\cos\theta,
\end{cases}\label{eq:pjn}
\end{equation}
the singular integral (\ref{eq:singInt}) leads back to:
\begin{align}
\varUpsilon_{sing}^{m,n} & =\mp\frac{1}{2\pi}\int_{0}^{l_{m}}dx\int_{-x\cos\theta}^{l_{n}-x\cos\theta}dx'\left[\frac{\cos\theta}{x'^{2}+x^{2}\sin^{2}\theta}-\frac{2x\sin^{2}\theta\left(x'+x\cos\theta\right)}{\left(x'^{2}+x^{2}\sin^{2}\theta\right)^{2}}\right]\nonumber \\
 & \times\left(1-\frac{x}{l_{m}}\right)\left[1-\frac{\left(x'+x\cos\theta\right)}{l_{n}}\right].\label{eq:singularAdjacentIntegral}
\end{align}
To carry out the inner integration, the following indefinite integrals
come in handy \cite{Gradshteyn2007}:
\begin{align}
\int\frac{dx'}{\left(x'^{2}+a^{2}\right)^{2}} & =\frac{1}{2a^{3}}\left[\frac{ax'}{x'^{2}+a^{2}}+\arctan\left(\frac{x'}{a}\right)\right];\\
\int\frac{x'dx'}{\left(x'^{2}+a^{2}\right)^{2}} & =-\frac{1}{2\left(x'^{2}+a^{2}\right)};\\
\int\frac{x'^{2}dx'}{\left(x'^{2}+a^{2}\right)^{2}} & =\frac{1}{2a}\left[\frac{-ax'}{x'^{2}+a^{2}}+\arctan\left(\frac{x'}{a}\right)\right];\\
\int\frac{x'^{3}dx'}{\left(x'^{2}+a^{2}\right)^{2}} & =\frac{1}{2}\left[\frac{a^{2}}{x'^{2}+a^{2}}+\log\left(x'^{2}+a^{2}\right)\right].
\end{align}
In particular, using \cite{Abramowitz1972}:
\begin{equation}
\arctan\left(z_{1}\right)\pm\arctan\left(z_{2}\right)=\arctan\left(\frac{z_{1}\pm z_{2}}{1\mp z_{1}z_{2}}\right),\label{eq:arctanSum}
\end{equation}
we obtain:
\begin{align}
\int_{-x\cos\theta}^{l_{n}-x\cos\theta}\frac{dx'}{\left(x'^{2}+x^{2}\sin^{2}\theta\right)^{2}} & =\frac{l_{n}\sin\theta\left(1-2\cos^{2}\theta\right)x+l_{n}^{2}\cos\theta\sin\theta}{2x^{3}\sin^{3}\theta\left(x^{2}-2l_{n}\cos\theta x+l_{n}^{2}\right)}+\frac{\arctan\left(\frac{l_{n}\sin\theta}{x-l_{n}\cos\theta}\right)}{2x^{3}\sin^{3}\theta};\\
\int_{-x\cos\theta}^{l_{n}-x\cos\theta}\frac{x'dx'}{\left(x'^{2}+x^{2}\sin^{2}\theta\right)^{2}} & =\frac{1}{2x^{2}}\left[\frac{-2l_{n}\cos\theta x+l_{n}^{2}}{x^{2}-2l_{n}\cos\theta x+l_{n}^{2}}\right];\\
\int_{-x\cos\theta}^{l_{n}-x\cos\theta}\frac{x'^{2}dx'}{\left(x'^{2}+x^{2}\sin^{2}\theta\right)^{2}} & =-\frac{l_{n}\sin\theta\left(1-2\cos^{2}\theta\right)x+l_{n}^{2}\cos\theta\sin\theta}{2x\sin\theta\left(x^{2}-2l_{n}\cos\theta x+l_{n}^{2}\right)}+\frac{\arctan\left(\frac{l_{n}\sin\theta}{x-l_{n}\cos\theta}\right)}{2x\sin\theta};\\
\int_{-x\cos\theta}^{l_{n}-x\cos\theta}\frac{x'^{3}dx'}{\left(x'^{2}+x^{2}\sin^{2}\theta\right)^{2}} & =\frac{1}{2}\left[\frac{2l_{n}\sin^{2}\theta\cos\theta x-l_{n}^{2}\sin^{2}\theta}{x^{2}-2l_{n}\cos\theta x+l_{n}^{2}}+\log\left(\frac{x^{2}-2l_{n}\cos\theta x+l_{n}^{2}}{x^{2}}\right)\right].
\end{align}
With the help of the previous formulas, equation (\ref{eq:singularAdjacentIntegral})
can be reduced to:
\begin{align}
\varUpsilon_{sing}^{m,n} & =\mp\frac{1}{2\pi}\left[1-\int_{\varepsilon}^{l_{m}}\frac{dx}{x}+\frac{\cos\theta}{2l_{n}}\int_{0}^{l_{m}}dx\log\left(\frac{x^{2}-2l_{n}\cos\theta x+l_{n}^{2}}{x^{2}}\right)\left(\frac{x}{l_{m}}-1\right)+\right.\nonumber \\
 & \left.-\frac{\sin\theta}{l_{n}}\int_{0}^{l_{m}}dx\arctan\left(\frac{l_{n}\sin\theta}{x-l_{n}\cos\theta}\right)\left(\frac{x}{l_{m}}-1\right)\right],\label{eq:simplifiedIntegral}
\end{align}
where the infinitesimal parameter $\varepsilon\rightarrow0$ has been
introduced in order to isolate the singularity arising from the integration
of the $-x^{-1}$ term, which turns out to be:
\begin{equation}
\mp\frac{1}{2\pi}\log\varepsilon.
\end{equation}
As expected, this singularity is just the opposite of that in (\ref{eq:upsilon11and22}).
Since the same result is obtained exchanging the elements $S_{n}$
and $S_{m}$, it is now clear that the divergences resulting from
the coincident integrations $\varUpsilon_{11}^{n}$, $\varUpsilon_{22}^{n}$
and from the adjacent integrations $\varUpsilon_{sing}^{m,n}$, $\varUpsilon_{sing}^{n,m}$
do indeed cancel in the sum (\ref{eq:hyperSingular_Matrix}) with
$i\equiv j$, and can therefore be explicitly removed. Finally, the
remaining terms in (\ref{eq:simplifiedIntegral}) can be expressed
analytically by means of standard integration formulas and algebraic
manipulation:%
\begin{align}
 & \varUpsilon_{sing}^{m,n}\pm\frac{1}{2\pi}\log\varepsilon=\mp\frac{1}{2\pi}\left\{ \frac{1}{2}+\frac{\sin\theta}{2l_{m}l_{n}}\left[l_{m}^{2}\arctan\left(\frac{l_{n}\sin\theta}{l_{m}-l_{n}\cos\theta}\right)+l_{n}^{2}\arctan\left(\frac{l_{m}\sin\theta}{l_{n}-l_{m}\cos\theta}\right)\right]+\right.\nonumber \\
 & \left.-\log\left(l_{m}l_{n}\right)-\frac{\cos\theta}{4l_{m}l_{n}}\left[\left(l_{m}^{2}+l_{n}^{2}-2l_{m}l_{n}\sec\theta\right)\log\left(l_{m}^{2}+l_{n}^{2}-2l_{m}l_{n}\cos\theta\right)-l_{m}^{2}\log\left(l_{m}^{2}\right)-l_{n}^{2}\log\left(l_{n}^{2}\right)\right]\right\} .\label{eq:analyticIntegral}
\end{align}

\section{Hypersingular operator: variational approach\label{sec:variationalApproach}}

In Section \ref{sec:directMethod}, a direct method to evaluate the
hypersingular integrals in (\ref{eq:hyperSingular_Matrix}) has been
proposed which makes use of an explicit cancellation of the residual
logarithmic divergences. Although such singularity subtraction procedure
does simplify the evaluation of the integrals over coincident elements,
rewritten through (\ref{eq:ups11and22regularized}) and (\ref{eq:upsilon12and21})
as a combination of the previously derived expressions (\ref{eq:I11I22})
and (\ref{eq:I12I21}) and of some well-known analytic functions,
the same does not hold for the case of adjacent segments, where it
requires the additional implementation of (\ref{eq:regInt}) and (\ref{eq:analyticIntegral}).
Due to this limitation, an alternative formulation based on the variational
approach described in \cite{Nedelec2001,Steinbach2008} will be considered
in the present section.

The hypersingular operator (\ref{eq:hyperSingular}) may be defined
more formally as the normal derivative of the double layer potential:
\begin{equation}
\hat{N}\left[f\right]\left(\mathbf{r}\right)\equiv\frac{\partial\hat{D}\left[f\right]\left(\mathbf{r}\right)}{\partial n}\equiv\lim_{\varepsilon\rightarrow0}\left\{ \mathbf{n}\cdot\nabla_{\mathbf{r}_{\varepsilon}}\left[\int_{S}d\mathbf{r}'\frac{\partial g\left(\mathbf{r}_{\varepsilon},\mathbf{r}'\right)}{\partial n'}f\left(\mathbf{r}'\right)\right]\right\} ,\label{eq:Nalternative}
\end{equation}
with: 
\begin{equation}
\mathbf{r}_{\varepsilon}\equiv\left(x_{\varepsilon},\,y_{\varepsilon}\right)\equiv\mathbf{r}+\varepsilon\mathbf{n}=\left(x+\varepsilon n_{x},\,y+\varepsilon n_{y}\right)
\end{equation}
representing a point in the tubular neighborhood of $S$. It is worth
noting that the appearance of divergent terms in the formulas of the
previous section could be interpreted as the effect of interchanging
the limit and the normal derivative in (\ref{eq:Nalternative}), as
the Cauchy principal value of the resulting integral is not defined.
Instead of focusing on (\ref{eq:Nalternative}), we may consider the
bilinear form induced by the hypersingular operator:
\begin{equation}
\left\langle \zeta\right|\hat{N}\left|\psi\right\rangle =\int_{S}d\mathbf{r}\,\zeta\left(\mathbf{r}\right)\hat{N}\left[\psi\right]\left(\mathbf{r}\right),\label{eq:bilinearForm}
\end{equation}
where $\psi\left(\mathbf{r}\right)$ and $\zeta\left(\mathbf{r}\right)$
are piecewise differentiable and globally continuous functions on
$S$. Exploiting the symmetry of the Green function $g\left(\mathbf{r},\mathbf{r}'\right)$
in (\ref{eq:Green}), it is easy to show that: 
\begin{align}
\frac{\partial}{\partial x}\left[\frac{\partial g\left(\mathbf{r},\mathbf{r}'\right)}{\partial n'}\right] & =-\mathbf{n}'\cdot\nabla'\left[\frac{\partial g\left(\mathbf{r},\mathbf{r}'\right)}{\partial x'}\right];\label{eq:xDerivative}\\
\frac{\partial}{\partial y}\left[\frac{\partial g\left(\mathbf{r},\mathbf{r}'\right)}{\partial n'}\right] & =-\mathbf{n}'\cdot\nabla'\left[\frac{\partial g\left(\mathbf{r},\mathbf{r}'\right)}{\partial y'}\right].\label{eq:yDerivative}
\end{align}
Now, introducing the following operator:
\begin{equation}
\mathrm{curl}_{S}\equiv\mathbf{n}\cdot\mathbf{curl}=n_{x}\frac{\partial}{\partial y}-n_{y}\frac{\partial}{\partial x},\label{eq:curlS}
\end{equation}
where $\mathbf{curl}$ is the surface curl on $\mathbb{R}^{2}$, and
making use of the Green function equation:
\begin{equation}
\Delta g\left(\mathbf{r},\mathbf{r}'\right)+k^{2}g\left(\mathbf{r},\mathbf{r}'\right)=\mp\delta\left(\mathbf{r}-\mathbf{r}'\right),
\end{equation}
we have:
\begin{align}
\mathrm{curl}_{S}'\left[\frac{\partial g\left(\mathbf{r}_{\varepsilon},\mathbf{r}'\right)}{\partial x'}\right] & =n_{x}'\frac{\partial^{2}g\left(\mathbf{r}_{\varepsilon},\mathbf{r}'\right)}{\partial y'\partial x'}-n_{y}'\frac{\partial^{2}g\left(\mathbf{r}_{\varepsilon},\mathbf{r}'\right)}{\partial x'^{2}}\nonumber \\
 & =n_{x}'\frac{\partial^{2}g\left(\mathbf{r}_{\varepsilon},\mathbf{r}'\right)}{\partial y'\partial x'}+n_{y}'\frac{\partial^{2}g\left(\mathbf{r}_{\varepsilon},\mathbf{r}'\right)}{\partial y'^{2}}-n_{y}'\Delta'g\left(\mathbf{r}_{\varepsilon},\mathbf{r}'\right)\nonumber \\
 & =n_{x}'\frac{\partial^{2}g\left(\mathbf{r}_{\varepsilon},\mathbf{r}'\right)}{\partial x'\partial y'}+n_{y}'\frac{\partial^{2}g\left(\mathbf{r}_{\varepsilon},\mathbf{r}'\right)}{\partial y'^{2}}+n_{y}'k^{2}g\left(\mathbf{r}_{\varepsilon},\mathbf{r}'\right)\nonumber \\
 & =\mathbf{n}'\cdot\nabla'\left[\frac{\partial g\left(\mathbf{r}_{\varepsilon},\mathbf{r}'\right)}{\partial y'}\right]+n_{y}'k^{2}g\left(\mathbf{r}_{\varepsilon},\mathbf{r}'\right).\label{eq:curlx}
\end{align}
Similarly:
\begin{equation}
\mathrm{curl}_{S}'\left[\frac{\partial g\left(\mathbf{r}_{\varepsilon},\mathbf{r}'\right)}{\partial y'}\right]=-\mathbf{n}'\cdot\nabla'\left[\frac{\partial g\left(\mathbf{r}_{\varepsilon},\mathbf{r}'\right)}{\partial x'}\right]-n_{x}'k^{2}g\left(\mathbf{r}_{\varepsilon},\mathbf{r}'\right).\label{eq:curly}
\end{equation}
Then, from (\ref{eq:xDerivative}) and (\ref{eq:curly}), we can write:
\begin{align}
\int_{S}d\mathbf{r}'\frac{\partial^{2}g\left(\mathbf{r}_{\varepsilon},\mathbf{r}'\right)}{\partial x_{\varepsilon}\partial n'}\psi\left(\mathbf{r}'\right) & =-\int_{S}d\mathbf{r}'\mathbf{n}'\cdot\nabla'\left[\frac{\partial g\left(\mathbf{r}_{\varepsilon},\mathbf{r}'\right)}{\partial x'}\right]\psi\left(\mathbf{r}'\right)\nonumber \\
 & =\int_{S}d\mathbf{r}'\mathrm{curl}_{S}'\left[\frac{\partial g\left(\mathbf{r}_{\varepsilon},\mathbf{r}'\right)}{\partial y'}\right]\psi\left(\mathbf{r}'\right)+k^{2}\int_{S}d\mathbf{r}'n_{x}'g\left(\mathbf{r}_{\varepsilon},\mathbf{r}'\right)\psi\left(\mathbf{r}'\right)\nonumber \\
 & =-\int_{S}d\mathbf{r}'\frac{\partial g\left(\mathbf{r}_{\varepsilon},\mathbf{r}'\right)}{\partial y'}\mathrm{curl}_{S}'\psi\left(\mathbf{r}'\right)+k^{2}\int_{S}d\mathbf{r}'n_{x}'g\left(\mathbf{r}_{\varepsilon},\mathbf{r}'\right)\psi\left(\mathbf{r}'\right),
\end{align}
where integration by parts has been applied in the last equality.
An analogous result is obtained from (\ref{eq:yDerivative}) and (\ref{eq:curlx}):
\begin{equation}
\int_{S}d\mathbf{r}'\frac{\partial^{2}g\left(\mathbf{r}_{\varepsilon},\mathbf{r}'\right)}{\partial y_{\varepsilon}\partial n'}\psi\left(\mathbf{r}'\right)=\int_{S}d\mathbf{r}'\frac{\partial g\left(\mathbf{r}_{\varepsilon},\mathbf{r}'\right)}{\partial x'}\mathrm{curl}_{S}'\psi\left(\mathbf{r}'\right)+k^{2}\int_{S}d\mathbf{r}'n_{y}'g\left(\mathbf{r}_{\varepsilon},\mathbf{r}'\right)\psi\left(\mathbf{r}'\right),
\end{equation}
so that, combining the two expressions under the limit $\varepsilon\rightarrow0$,
we get:
\begin{align*}
\hat{N}\left[\psi\right]\left(\mathbf{r}\right) & =\fint_{S}d\mathbf{r}'\mathrm{curl}_{S}g\left(\mathbf{r},\mathbf{r}'\right)\mathrm{curl}_{S}'\psi\left(\mathbf{r}'\right)+k^{2}\fint_{S}d\mathbf{r}'g\left(\mathbf{r},\mathbf{r}'\right)\psi\left(\mathbf{r}'\right)\left(\mathbf{n}\cdot\mathbf{n}'\right).
\end{align*}
Finally, using again integration by parts, (\ref{eq:bilinearForm})
is reduced to:
\begin{equation}
\left\langle \zeta\right|\hat{N}\left|\psi\right\rangle =\int_{S}d\mathbf{r}\fint_{S}d\mathbf{r}'g\left(\mathbf{r},\mathbf{r}'\right)\left[k^{2}\left(\mathbf{n}\cdot\mathbf{n}'\right)\zeta\left(\mathbf{r}\right)\psi\left(\mathbf{r}'\right)-\mathrm{curl}_{S}\zeta\left(\mathbf{r}\right)\mathrm{curl}_{S}'\psi\left(\mathbf{r}'\right)\right].\label{eq:solvedBilinearForm}
\end{equation}
In other words, the bilinear form induced by the hypersingular operator
has been recast as a bilinear form induced by the single layer potential. 

It is apparent that the matrix $N_{ij}$ defined in (\ref{eq:hyperSingular_Matrix})
may be seen as a discrete version of the bilinear form (\ref{eq:bilinearForm})
where $\zeta\left(\mathbf{r}\right)$ and $\psi\left(\mathbf{r}'\right)$
are replaced by linear or higher-order basis functions $p_{i}\left(\mathbf{r}\right)$
and $p_{j}\left(\mathbf{r}'\right)$, respectively. Therefore,
from (\ref{eq:solvedBilinearForm}), it follows that:
\begin{align}
N_{ij} & \equiv\sum_{m\in i}\,\sum_{n\in j}\int_{S_{m}}d\mathbf{r}\fint_{S_{n}}d\mathbf{r}'g\left(\mathbf{r},\mathbf{r}'\right)\left[k^{2}\left(\mathbf{n}\cdot\mathbf{n}'\right)p_{i}^{m}\left(\mathbf{r}\right)p_{j}^{n}\left(\mathbf{r}'\right)-\mathrm{curl}_{S_{m}}p_{i}^{m}\left(\mathbf{r}\right)\mathrm{curl}_{S_{n}}'p_{j}^{n}\left(\mathbf{r}'\right)\right].\label{eq:solvedHypersingularMatrix}
\end{align}
In order to apply the $\mathrm{curl}_{S}$ operator to the triangular
functions (\ref{eq:basisFunctions}), we start by lifting the parameterization
(\ref{eq:parametricCoordinates}) into the two-dimensional tubular
neighborhood of the $n$-th segment:
\begin{equation}
\mathbf{r}_{\varepsilon}\left(t_{n}\right)=\mathbf{r}_{A}^{n}+\left(\mathbf{r}_{B}^{n}-\mathbf{r}_{A}^{n}\right)t_{n}+\varepsilon\mathbf{n}.\label{eq:parametricExtension}
\end{equation}
Then, taking the inner product of (\ref{eq:parametricExtension})
with $\mathbf{l}_{n}\equiv\left(\mathbf{r}_{B}^{n}-\mathbf{r}_{A}^{n}\right)$:
\begin{equation}
\mathbf{l}_{n}\cdot\left[\mathbf{r}_{\varepsilon}\left(t_{n}\right)-\mathbf{r}_{A}^{n}\right]=l_{n}^{2}t_{n},
\end{equation}
we can write:
\begin{equation}
p_{j}^{n}\left(x,y\right)\equiv\begin{cases}
1-\frac{1}{l_{n}^{2}}\left[\left(x-x_{A}^{n}\right)\left(x_{B}^{n}-x_{A}^{n}\right)+\left(y-y_{A}^{n}\right)\left(y_{B}^{n}-y_{A}^{n}\right)\right] & \mathrm{if}\:\mathbf{r}_{j}=\mathbf{r}_{A}^{n};\\
\frac{1}{l_{n}^{2}}\left[\left(x-x_{A}^{n}\right)\left(x_{B}^{n}-x_{A}^{n}\right)+\left(y-y_{A}^{n}\right)\left(y_{B}^{n}-y_{A}^{n}\right)\right] & \mathrm{if}\:\mathbf{r}_{j}=\mathbf{r}_{B}^{n}.
\end{cases}\label{eq:extendedBasisFunctions}
\end{equation}
Equation (\ref{eq:extendedBasisFunctions}) provides a constant extension
of the functions (\ref{eq:basisFunctions}) along $\mathbf{n}$.
On using (\ref{eq:curlS}) and the definition of unit normal to the
$n$-th segment in $\mathbb{R}^{3}$:
\begin{equation}
\mathbf{n}=\frac{\mathbf{l}_{n}\times\mathbf{z}}{\left|\mathbf{l}_{n}\times\mathbf{z}\right|}=\frac{1}{l_{n}}\left(y_{B}^{n}-y_{A}^{n},\,x_{A}^{n}-x_{B}^{n},\,0\right),
\end{equation}
it follows that:
\begin{equation}
\mathrm{curl}_{S_{n}}p_{j}^{n}\left(x,y\right)=\begin{cases}
-\frac{1}{l_{n}} & \mathrm{if}\:\mathbf{r}_{j}=\mathbf{r}_{A}^{n};\\
\frac{1}{l_{n}} & \mathrm{if}\:\mathbf{r}_{j}=\mathbf{r}_{B}^{n}.
\end{cases}
\end{equation}

All the integrations in (\ref{eq:solvedHypersingularMatrix}) can
now be computed as in Section \ref{sec:singleLayer}. In particular,
Gauss-Legendre quadrature formulas directly apply whenever $S_{m}\neq S_{n}$.
On the other hand, making use of (\ref{eq:Ialphan}) and (\ref{eq:IalphaInt}),
the four possible integrals over coincident segments acquire the following
form:
\begin{align}
\widetilde{\varUpsilon}_{11}^{n} & =\widetilde{\varUpsilon}_{22}^{n}=k^{2}I_{11}^{n}-\frac{1}{l_{n}^{2}}\int_{S_{n}}d\mathbf{r}\fint_{S_{n}}d\mathbf{r}'g\left(\mathbf{r},\mathbf{r}'\right)=k^{2}I_{11}^{n}-\frac{i}{2k^{2}l_{n}^{2}}\varGamma_{0}\left(kl_{n}\right);\label{eq:upsilon1122tilda}\\
\widetilde{\varUpsilon}_{12}^{n} & =\widetilde{\varUpsilon}_{21}^{n}=k^{2}I_{12}^{n}+\frac{1}{l_{n}^{2}}\int_{S_{n}}d\mathbf{r}\fint_{S_{n}}d\mathbf{r}'g\left(\mathbf{r},\mathbf{r}'\right)=k^{2}I_{12}^{n}+\frac{i}{2k^{2}l_{n}^{2}}\varGamma_{0}\left(kl_{n}\right),\label{eq:upsilon1221tilda}
\end{align}
where $I_{11}^{n}$, $I_{12}^{n}$ and $\varGamma_{0}\left(kl_{n}\right)$
are defined in (\ref{eq:I11I22})-(\ref{eq:gamma0}) and a tilde has
been introduced to avoid notation overlap. It is important to understand
that, despite the similarity between the last two expressions and
(\ref{eq:upsilon11and22})-(\ref{eq:upsilon12and21}), a comparison
of the direct and variational methods is only possible in terms of
the matrix entries $N_{ij}$, that is to say, after summing all the
four integrals over $\left\{ m\in i\right\} \times\left\{ n\in j\right\} $.
An example in this regard is reported in Table \ref{tab:comparison}.
\begin{table}
\begin{centering}
\begin{tabular}{|c|c|c|c|}
\hline 
$k$ & %
\begin{tabular}{c}
\noalign{\vskip\doublerulesep}
coincident nodes $i$, $j$\tabularnewline[\doublerulesep]
\end{tabular} & %
\begin{tabular}{c}
\noalign{\vskip\doublerulesep}
first-neighbor nodes $i$, $j$\tabularnewline[\doublerulesep]
\end{tabular} & %
\begin{tabular}{c}
\noalign{\vskip\doublerulesep}
second-neighbor nodes $i$, $j$\tabularnewline[\doublerulesep]
\end{tabular}\tabularnewline
\hline 
{\footnotesize 0.1} & {\footnotesize%
\begin{tabular}{c}
\noalign{\vskip0.3cm}
$N_{ij}^{dir}=-0.440599+0.004556\,i$\tabularnewline[0.15cm]
\noalign{\vskip0.15cm}
$N_{ij}^{var}=-0.440600+0.004556\,i$\tabularnewline[0.3cm]
\end{tabular}} & {\footnotesize%
\begin{tabular}{c}
\noalign{\vskip0.3cm}
$N_{ij}^{dir}=0.056506+0.004426\,i$\tabularnewline[0.15cm]
\noalign{\vskip0.15cm}
$N_{ij}^{var}=0.056461+0.004426\,i$\tabularnewline[0.3cm]
\end{tabular}} & {\footnotesize%
\begin{tabular}{c}
\noalign{\vskip0.3cm}
$N_{ij}^{dir}=0.065692+0.005318\,i$\tabularnewline[0.15cm]
\noalign{\vskip0.15cm}
$N_{ij}^{var}=0.065691+0.005318\,i$\tabularnewline[0.3cm]
\end{tabular}}\tabularnewline
\hline 
{\footnotesize 1} & {\footnotesize%
\begin{tabular}{c}
\noalign{\vskip0.3cm}
$N_{ij}^{dir}=-0.155208+0.396113\,i$\tabularnewline[0.15cm]
\noalign{\vskip0.15cm}
$N_{ij}^{var}=-0.155211+0.396113\,i$\tabularnewline[0.3cm]
\end{tabular}} & {\footnotesize%
\begin{tabular}{c}
\noalign{\vskip0.3cm}
$N_{ij}^{dir}=0.111035+0.277807\,i$\tabularnewline[0.15cm]
\noalign{\vskip0.15cm}
$N_{ij}^{var}=0.110990+0.277807\,i$\tabularnewline[0.3cm]
\end{tabular}} & {\footnotesize%
\begin{tabular}{c}
\noalign{\vskip0.3cm}
$N_{ij}^{dir}=-0.037765+0.011613\,i$\tabularnewline[0.15cm]
\noalign{\vskip0.15cm}
$N_{ij}^{var}=-0.037766+0.011613\,i$\tabularnewline[0.3cm]
\end{tabular}}\tabularnewline
\hline 
{\footnotesize 10} & {\footnotesize%
\begin{tabular}{c}
\noalign{\vskip0.3cm}
$N_{ij}^{dir}=0.014973+6.447260\,i$\tabularnewline[0.15cm]
\noalign{\vskip0.15cm}
$N_{ij}^{var}=0.014811+6.447260\,i$\tabularnewline[0.3cm]
\end{tabular}} & {\footnotesize%
\begin{tabular}{c}
\noalign{\vskip0.3cm}
$N_{ij}^{dir}=-0.008735+2.135620\,i$\tabularnewline[0.15cm]
\noalign{\vskip0.15cm}
$N_{ij}^{var}=-0.008780+2.135620\,i$\tabularnewline[0.3cm]
\end{tabular}} & {\footnotesize%
\begin{tabular}{c}
\noalign{\vskip0.3cm}
$N_{ij}^{dir}=-0.001017+0.013765\,i$\tabularnewline[0.15cm]
\noalign{\vskip0.15cm}
$N_{ij}^{var}=-0.001017+0.013765\,i$\tabularnewline[0.3cm]
\end{tabular}}\tabularnewline
\hline 
\end{tabular}
\par\end{centering}
\caption{Numerical comparison between the direct method presented in Section
\ref{sec:directMethod} and the variational approach of Section \ref{sec:variationalApproach}
by evaluation of an arbitrary hypersingular matrix entry $N_{ij}$ over a
mesh like that shown in figures \ref{fig:segments} and \ref{fig:combinations}
(average length of the segments $\approx2.26$). The number of integration
points has been set to $20$ for all Gauss-Legendre quadratures. \label{tab:comparison}}

\end{table}

\section{Conclusions\label{sec:conclusions}}

In this work, extensive use of analytic integration has been made
to provide quasi-closed-form expressions for the Galerkin singular
integrals of the Helmholtz boundary operators in two dimensions. Two
different techniques have been applied to the discrete hypersingular
operator, namely, a direct method and a variational formulation; the
second approach proves superior in that it does not require singularity
subtraction. To summarize, the relevant formulas are given by (\ref{eq:I11I22})-(\ref{eq:I12I21})
for the single layer operator and by (\ref{eq:upsilon1122tilda})-(\ref{eq:upsilon1221tilda})
for the hypersingular operator, and they rely on the numerical evaluation
of well-known analytic functions and of the integrals (\ref{eq:gamma0})-(\ref{eq:gamma2}).
These formulas may simplify the implementation of the BEM in two-dimensional
electromagnetic, acoustic and quantum mechanical problems.

\ack
I would like to acknowledge Prof. Francesco Andriulli, the members of the Computational 
Electromagnetics Research Laboratory and all my colleagues from the Microwaves Department 
of IMT Atlantique, for their support and friendship.

\appendix

\section{Singularity cancellation for the double layer adjacent integrations\label{sec:appendix}}

By choosing $\mathbf{r}_{0}$ to represent the position of the
common vertex and $\mathbf{l}_{m}$ and $\mathbf{l}_{n}$
the distance vectors between $\mathbf{r}_{0}$ and the outer extrema
of $S_{m}$ and $S_{n}$, respectively, we can define the local variables
$t,\:t'\in\left[0,1\right]$ such that:
\begin{equation}
\mathbf{r}\left(t\right)=\mathbf{r}_{0}+\mathbf{l}_{m}t;\qquad\mathbf{r}'\left(t'\right)=\mathbf{r}_{0}+\mathbf{l}_{n}t';\qquad\mathbf{R}\left(t,t'\right)=\mathbf{l}_{m}t-\mathbf{l}_{n}t'.
\end{equation}
Now, applying the coordinate transformation: 
\begin{equation}
\begin{cases}
t=\rho\cos\phi;\\
t'=\rho\sin\phi,
\end{cases}
\end{equation}
the singular double integral over adjacent segments in (\ref{eq:doubleLayer_Matrix})
can be rewritten as:
\begin{equation}
\left[\int_{0}^{\pi/4}d\phi\int_{0}^{\sec\phi}d\rho+\int_{\pi/4}^{\pi/2}d\phi\int_{0}^{\csc\phi}d\rho\right]\:f_{D}\left(\rho,\phi\right),\label{eq:integral}
\end{equation}
where the integrand function $f_{D}\left(\rho,\phi\right)$ is defined
by:
\begin{align}
f_{D}\left(\rho,\phi\right) & =-\frac{ikl_{m}^{2}l_{n}\sin\theta}{4}\left\{ \frac{H_{1}^{(1,2)}\left[kR\left(\rho,\phi\right)\right]}{R\left(\rho,\phi\right)}\rho^{2}\cos\phi\left(1-\rho\cos\phi\right)\left(1-\rho\sin\phi\right)\right\} ,\label{eq:integrand}
\end{align}
and the following expressions have been considered:
\begin{align}
 & \mathbf{l}_{m}\cdot\mathbf{l}_{n}=l_{n}l_{m}\cos\theta;\qquad\mathbf{R}\cdot\mathbf{n}'=\mathbf{l}_{m}\cdot\mathbf{n}'\rho\cos\phi=-l_{m}\sin\theta\rho\cos\phi;\\
 & R\left(\rho,\phi\right)=\left|\mathbf{R}\right|=\rho\sqrt{l_{m}^{2}\cos^{2}\phi+l_{n}^{2}\sin^{2}\phi-2l_{n}l_{m}\cos\theta\cos\phi\sin\phi},
\end{align}
with $\mathbf{n}'$ oriented as in Figure \ref{fig:referenceFrames}
(right side).

Owing to the presence of the multiplicative factor $\rho$ from the
Jacobian, the integrand (\ref{eq:integrand}) is now regular at the
common vertex (namely, at $\rho=0$). However, the integration domain
is no longer rectangular in the new coordinates, so that both integrals
in (\ref{eq:integral}) need to be further transformed in order for
Gauss-Legendre quadrature to apply. The required variable changes
are easily shown to be $\rho\rightarrow\rho\cos\phi$ for the first
integral, and $\rho\rightarrow\rho\sin\phi$ for the second, leading
to:
\begin{equation}
\int_{0}^{\pi/4}d\phi\int_{0}^{1}d\rho\,f_{D}\left(\rho\sec\phi,\phi\right)\sec\phi+\int_{\pi/4}^{\pi/2}d\phi\int_{0}^{1}d\rho\,f_{D}\left(\rho\csc\phi,\phi\right)\csc\phi,
\end{equation}
both of which can now be solved numerically.\\\\


\providecommand{\newblock}{}

\end{document}